# ON THE STRUCTURE OF QUASI-STATIONARY COMPETING PARTICLE SYSTEMS[1]

By Louis-Pierre Arguin and Michael Aizenman

*Princeton University*

We study point processes on the real line whose configurations $X$ are locally finite, have a maximum and evolve through increments which are functions of correlated Gaussian variables. The correlations are intrinsic to the points and quantified by a matrix $Q = \{q_{ij}\}_{i,j \in \mathbb{N}}$. A probability measure on the pair $(X, Q)$ is said to be *quasi-stationary* if the joint law of the gaps of $X$ and of $Q$ is invariant under the evolution. A known class of universally quasi-stationary processes is given by the Ruelle Probability Cascades (RPC), which are based on hierarchically nested Poisson–Dirichlet processes. It was conjectured that up to some natural superpositions these processes exhausted the class of laws which are robustly quasi-stationary. The main result of this work is a proof of this conjecture for the case where $q_{ij}$ assume only a finite number of values. The result is of relevance for mean-field spin glass models, where the evolution corresponds to the cavity dynamics, and where the hierarchical organization of the Gibbs measure was first proposed as an ansatz.

## 1. Introduction.

1.1. *The problem.* A competing particle system is a point process on $\mathbb{R}$ whose configurations can be ordered, that is have a maximum and are locally finite, and which evolve here in discrete time steps. We are interested in the situation where the points represent the location of entities which have intrinsic characteristics that are of relevance for their dynamics. In the case discussed here, these are not affected by the time evolution and in particular by the particles' evolving positions on the line. The evolution is random, with steps which in any given time-increment are correlated through the

Received September 2007; revised September 2008.

[1]Supported in part by NSERC and FQRNT postgraduate fellowships, and NSF Grant DMS-06-02360.

points' characteristics. The process is said to be *quasi-stationary* if the joint distribution of the gaps among the points along with their characteristic features is stationary in time. Our goal is to characterize processes which are endowed with a robust quasi-stationarity property. We prove, under some simplifying restrictions, that such a condition requires the process to be organized as a hierarchical random probability cascade, RPC—an acronym which can also be read as Ruelle Probability Cascade. This remarkable family of processes was introduced in the study of spin-glass models [9, 15, 16], for which it plays a fundamental role [1], and it has attracted a great deal of attention since [8].

More explicitly, in the system under consideration, the configuration of the process is given by:

1. an ordered sequence of positions, $X = \{X_i\}_{i \in \mathbb{N}} \subset \mathbb{R}$, with $X_1 \geq X_2 \geq \cdots$, and

2. a so-called *overlap* matrix, which is a positive definite quadratic form given by $Q = \{q_{ij}\}$, with $q_{ii} = 1$ for all $i \in \mathbb{N}$ (and hence $|q_{ij}| \leq 1$).

The time evolution is given by correlated increments, $X_i \mapsto X_i + \psi(\kappa_i)$, where $\kappa$ is a Gaussian field on $\mathbb{N}$, independent of $X$ with covariance given by the matrix $Q$, and $\psi$ is a real function. This is followed by a reordering of the indices, denoted by $\downarrow$, so that the full evolution step is

$$(1.1) \qquad\qquad (X, Q) \mapsto (X', Q')$$

with

$$
\begin{aligned}
(1.2) \quad & X' = (X_i + \psi(\kappa_i), i \in \mathbb{N})_\downarrow = (X_{\pi^{-1}(m)} + \psi(\kappa_{\pi^{-1}(m)}), m \in \mathbb{N}), \\
& Q' = \{q_{\pi^{-1}(m)\pi^{-1}(n)}\},
\end{aligned}
$$

where $\pi$ is the permutation of $\mathbb{N}$ induced by $\downarrow$: $\pi(i) = m$ if the $i$th point in $X$ is mapped to the $m$th point in $X'$. If such a permutation does not exist, the evolution is not defined.

It may be remarked here that with a suitable choice of $\psi$ any random variable may be presented as a function of a Gaussian, $\psi(\kappa)$. This formulation is, however, convenient for the description of the structure of correlations between the variables encountered here. With the assumptions which are spelled below on $\psi$ it also carries some mild implications on the distribution of a single incremental variable (continuity of its distribution, and finiteness of its moment generating function).

The above evolution can be regarded as a competition within a crowd of points. While the labels (i) change according to the relative positions, the covariances $q_{ij}$ are intrinsic to the particles. A crucial feature of this dynamical system is that both the position and the matrix of overlaps are random—the relevant stochastic object is the pair $(X, Q)$.

The present work addresses the question of describing the collection of distributions of $(X, Q)$ which are *quasi-stationary* under the above evolution, in the sense that the joint law of the gaps, $\{X_i - X_{i+1}\}_{i \in \mathbb{N}}$, and $q_{ij}$ is invariant under the evolution. In the particular case where $Q$ is the identity, that is the increments are independent, Ruzmaikina and Aizenman [17] have shown that a quasi-stationary process $X$ must be a superposition of Poisson processes with exponential density. The aforementioned Ruelle probability cascades form a broader class of processes, with $Q$ not limited to the identity matrix. Each process in this collection is quasi-stationary under any dynamics of the above form with an exponentially bounded $\psi$, and the class is parametrized by probability measures on $[0, 1]$. It was conjectured that this collection includes all the processes which are robustly quasi-stationary. The significance of this conjecture for spin glass models is discussed in [2].

The RPC measures incorporate an interesting hierarchical structure. That can alternatively be expressed by the statement that the metric induced by the overlaps, $d_{ij} = \sqrt{1 - q_{ij}}$, is almost surely ultrametric, in the sense that

$$d_{ij} \leq \max\{d_{ik}, d_{kj}\}$$

for any triplet $i, j, k$. The inequality can equivalently be expressed by saying that the condition $q_{ij} \geq r$ is transitive, as a condition on pairs $\{i, j\} \subset \mathbb{N}$. In such case, the set of points endowed with the distance $d$ has a tree structure and the overlap matrix is the covariance matrix of the Gaussian field on the tree.

The main result of this work establishes the above conjecture in the case of systems that are quasi-stationary in a robust sense and for which the set of values taken by the covariances of the point $i$ is finite. The result is stated precisely in Theorem 1.8, in the next section, where the necessary concepts are introduced.

Our interest in this problem is partially motivated by the desire to shed light on the properties of mean-field spin glass systems. It is believed that the equilibrium states of a wide class of such models are organized in an ultrametric fashion (among them the Sherrington–Kirkpatrick model for which this assumption leads to the correct free energy per particle [2, 12, 19]). It has been an interesting question of what is at the root of that. Since it makes sense to expect the law of the corresponding Gibbs states to be asymptotically invariant under a so-called *cavity dynamics* (see, e.g., [2, 11]), which has exactly the form (1.4), the result presented here may provide a step toward the explanation of this phenomenon.

**1.2.** *Competing particle systems and random overlap structures.* A convenient setup for our discussion is to regard the random pair $(X, Q)$ as a *random overlap structure* or *ROSt*. This concept was introduced in the study of mean-field spin glass systems where it provides a useful framework [1]. Its significance for the Parisi solution of the SK model was discussed in [2, 4].

The requirement that the evolving configuration of particles will admit sequential order with probability one is assured by the convenient condition, which is also natural for our result, that for some $\beta > 0$ the point process satisfies

$$(1.3) \qquad \sum_{i \in \mathbb{N}} e^{\beta X_i} < \infty \qquad \text{a.s.}$$

(For finite configurations, it should be understood that $X_i = -\infty$ for all but a finite number of points.)

Point processes of the above kind are closely related to random mass-partitions [7]. A *mass-partition* is a sequence $(s_i, i \in \mathbb{N})$ such that $s_1 \geq s_2 \geq \cdots \geq 0$ and $\sum_i s_i \leq 1$. It is said to be *proper* if $\sum_i s_i = 1$. The space $P_m$ of mass-partition is usually endowed with the metric $d(s, s') = \max_{i \in \mathbb{N}} \{|s_i - s_i'|\}$. Under this metric, $P_m$ is a compact separable space. The matrix $Q$ associated to $X$ takes value in the space $\mathcal{Q}$ of real symmetric positive semi-definite forms $\mathbb{N} \times \mathbb{N}$ for which the diagonal is normalized to 1. Elements of this space are called *overlap matrices* by analogy with spin glass systems. We endow this space with the topology it inherits as a subset of the compact space $[-1, 1]^{\mathbb{N} \times \mathbb{N}}$ equipped with the product topology. This space is then metrizable, compact and separable.

DEFINITION 1.1 (*ROSt* [2]). The space of overlap structures $\Omega_{os}$ is the compact separable metric space $P_m \times \mathcal{Q}$. We write $\mathcal{F}_{os}$ for the Borel $\sigma$-algebra on $\Omega_{os}$. A random overlap structure or ROSt is an element of $\mathcal{M}_1(\Omega_{os})$, the space of regular probability measure on $(\Omega_{os}, \mathcal{F}_{os})$, for which $Q$ is almost surely positive definite.

We will usually write $\mathbb{P}$ for the law of a random overlap structure $(\xi, Q)$ and $\mathbb{E}$ for its expectation. A basis for the topology of $\Omega_{os}$ are the open sets of the form

$$B(n, \{I_i\}, \{A_{ij}\}) = \{(s, Q) : s_i \in I_i, q_{ij} \in A_{ij} \ \forall 1 \leq i, j \leq n\}$$

for some $n \in \mathbb{N}$ and open sets $I_i \subset [0, 1]$, $A_{ij} \subset [-1, 1]$. It is clear from the basis generating $\mathcal{F}_{os}$ that the law $\mathbb{P}$ of a ROSt is determined by the expected value of measurable functions that depend only on a finite number of points. This class can be further restricted to continuous functions by a simple application of the Stone–Weierstrass theorem, and one has:

---

[2]The notion of random overlap structures [1, 2] could conceivably be discussed in the more general case, where $\xi$ is taken to be a random finite measure on a measure space $\mathcal{A}$ and $Q$ a positive semi-definite form on $\mathcal{A}$. Our definition corresponds to the case $\mathcal{A} = \mathbb{N}$.

PROPOSITION 1.2. *Let $(\xi, Q)$ and $(\xi', Q')$ be two ROSt's. Let $C_n(\Omega_{os})$ be the class of continuous functions on $\Omega_{os}$ that only depend on the positions and the overlaps of the first $n$ points. If $\mathbb{E}[f(\xi, Q)] = \mathbb{E}[f(\xi', Q')]$ for every $f \in C_n(\Omega_{os})$ and $n \in \mathbb{N}$, then $(\xi, Q)$ and $(\xi', Q')$ have the same law.*

To every pair $(X, Q)$ with $X$ satisfying the condition (1.3), we associate a ROSt $(\xi, Q)$ by simply taking, for some $\beta > 0$:

$$\xi_i := \frac{e^{\beta X_i}}{\sum_j e^{\beta X_j}}.$$

Because of the normalization, the law of $\xi$ depends only on the law of the gaps of $X$. Hence quasi-stationarity of $(X, Q)$ under the correlated evolution (1.1) is equivalent to the invariance of the law of $(\xi, Q)$ under the stochastic map

$$(1.4) \qquad (\xi, Q) \mapsto \Phi_{\psi(\kappa)}(\xi, Q) := \left( \left( \frac{\xi_i e^{\psi(\kappa_i)}}{\sum_j \xi_j e^{\psi(\kappa_j)}}, i \in \mathbb{N} \right)_{\downarrow}, \pi Q \pi^{-1} \right),$$

where, for simplicity, we incorporated the factor $\beta$ in $\psi$. The symbol $\downarrow$ means that the weights are reordered in decreasing order after evolution. Again, the permutation $\pi$ reflects the reordering, namely $\pi(i) = m$ if the $i$th weight becomes the $m$th weight after evolution. From now on, we will identify a permutation matrix, $\pi$, with the corresponding permutation matrix, $\pi_{ij} = 1$ if $\pi(i) = j$ and 0 otherwise. With this notation, the evolved overlap matrix $\{q_{\pi^{-1}(m)\pi^{-1}(n)}\}$ becomes $\pi Q \pi^{-1}$. A sufficient condition for the mapping to be nonsingular, that is, with $\sum_i \xi_i e^{\psi(\kappa_i)}$ finite a.s., is the finiteness of the expectation of $e^{\psi(\kappa_i)}$. In fact, we will assume the following condition, which suffices for this purpose, and also guarantees that the increments admit a continuous range of values (and in particular are not supported on a lattice).

ASSUMPTIONS 1.3. *The function $\psi : \mathbb{R} \to \mathbb{R}$ is a Borel measurable function satisfying*

$$\int_{\mathbb{R}} \frac{e^{-z^2/2}}{\sqrt{2\pi}} e^{\lambda \psi(z)} \, dz < \infty$$

*for any $\lambda \in \mathbb{R}$. Furthermore, for $Y$ a standard Gaussian variable the law of $\psi(Y)$ is absolutely continuous with respect to the Lebesgue measure.*

Assumption 1.3 is in particular fulfilled by any function in $C^2(\mathbb{R})$ with bounded derivatives, for example: $\psi(\kappa) = \kappa$ and $\psi(\kappa) = \log \cosh \kappa$. In the linear case, the evolution of $\log \xi$ is by Gaussian increments. The second example is of particular interest for the SK spin glass model.

It is straightforward to check that the mapping $\Phi_{(\cdot)}(\cdot)$ in (1.4) as a mapping from $\mathbb{R}^{\mathbb{N}} \times \Omega_{os}$ to $\Omega_{os}$ is jointly measurable. This guarantees that $\Phi_{\psi(\kappa)}(\xi, Q)$ is well defined as a ROSt induced by the laws of $\kappa$ and $(\xi, Q)$.

There are natural variations of the mapping (1.4) that will be of importance in our study. Let $Q^{*r}$ denote the $r$th power, $r \in \mathbb{N}$, of the matrix $Q$ in the sense of the entry-wise product. By a known theorem of Schur, if $Q$ is positive definite so is $Q^{*r}$ [14]. The stochastic mapping (1.4) can then be considered for a Gaussian field $\kappa$ with covariance $Q^{*r}$ for some $r \in \mathbb{N}$.

DEFINITION 1.4 (*Correlated evolution*). A correlated evolution on $\Omega_{os}$ is a stochastic mapping of the form (1.4), for which the increments are $\psi(\kappa)$ with $\kappa$ a centered Gaussian field with covariance $\mathbb{E}[\kappa_i \kappa_j] = g(q_{ij})$, for some $g \in \mathcal{C}(\mathbb{R})$ and $\psi$ a function satisfying Assumption 1.3.

As the function $\psi$ will often be fixed, we will sometimes drop the dependence on $\psi$ and the Gaussian field in the notation. We will write $\Phi_r$ for the correlated evolution with $g(q) \equiv q^r$, corresponding to the covariance $Q^{*r}$. Moreover, $\mathbb{P}_r$ will denote the probability measure on $\Omega_{os} \times \mathbb{R}^{\mathbb{N}}$ given by

(1.5)               $$d\mathbb{P}_r(\xi, Q, \kappa) := d\mathbb{P}(\xi, Q)\, d\nu_{Q^{*r}}(\kappa),$$

where $\nu_{Q^{*r}}$ is the Gaussian measure with covariance $Q^{*r}$.

Our original challenge translates into characterizing the random overlap structures that are quasi-stationary under some correlated evolution.

DEFINITION 1.5 (*Quasi-stationarity*). Let $\psi$ be a function satisfying Assumption 1.3:

- A ROSt $(\xi, Q)$ is said to be quasi-stationary under the correlated evolution $\Phi_r$ if and only if

$$\Phi_r(\xi, Q) \stackrel{\mathcal{D}}{=} (\xi, Q),$$

  where the symbol $\stackrel{\mathcal{D}}{=}$ means equality in distribution.
- A ROSt $(\xi, Q)$ is said to be robustly quasi-stationary if and only if it is quasi-stationary under $\Phi_r$ for $r = 1$ and an infinite number of $r \in \mathbb{N}$.

The ROSt consisting of only one point is trivially quasi-stationary. Furthermore, it is easy to see that under the dynamics discussed above any point process with finitely many points will converge to this trivial case (the distribution of the gaps would spread indefinitely). Hence, we focus our attention on systems with infinitely many particles.

For a fixed $\psi$ the subset of robustly quasi-stationary laws of $\mathcal{M}_1(\Omega_{os})$ is plainly convex. It can be checked that it is also closed in the weak topology induced by $C(\Omega_{os})$ on $\mathcal{M}_1(\Omega_{os})$. As $\Omega_{os}$ is a separable compact Hausdorff space, Choquet's theorem ensures that any robustly quasi-stationary law is a linear superposition of the extreme measures (see, e.g., [18], page 63). In our case, these extreme or *ergodic* measures are exactly the ones for which the $\Phi_r$-invariant functions are constant for all $r \in \mathbb{N}$. The $\Phi_r$-invariant functions are the bounded measurable functions $f : \Omega_{os} \to \mathbb{R}$ such that

$$\mathbb{E}_r[f(\Phi_r(\xi, Q))|(\xi, Q)] \overset{\mathbb{P}\text{-a.s.}}{=} f(\xi, Q),$$

where $\mathbb{E}[\cdot|X]$ denotes the conditional expectation with respect to the $\sigma$-algebra generated by the random variable $X$.

1.3. *Q-factorization.* For a given ROSt $(\xi, Q)$, we denote

$$S_Q := \{q_{ij} : 1 \le i < j < \infty\}$$

and

$$S_Q(i) := \{q_{ij} : j \ne i\}.$$

A ROSt is said to have a *finite state space* if $|S_Q| < \infty$ a.s. In this work we consider only such systems. Note that the cardinality of the state space $S_Q$ is a $\Phi_r$-invariant function for all $r \in \mathbb{N}$. If $(\xi, Q)$ is ergodic then $S_Q$ is a deterministic set.

DEFINITION 1.6. A ROSt $(\xi, Q)$ is said to be *overlap-indecomposable*, or simply *indecomposable*, if

$$S_Q(i) = S_Q \qquad \text{for all } i \in \mathbb{N}.$$

The set $S_Q(i)$ provides a tag according to which the points may be partitioned in a time independent fashion. Precisely, for any subset $A$ of $S_Q$, one defines $I_A := \{i : S_Q(i) = A\}$. A nonempty element $I_A$ in this partition gives rise to a ROSt as follows

$$(\xi, Q)_A := ((\xi_i, i \in I_A), \{q_{ij}\}_{i,j \in I_A}).$$

Each of these would be quasi-stationary if $(\xi, Q)$ was. ROSt's $(\xi, Q)_A$ need not be indecomposable, since the collection of overlaps of $i$ within the limited collection $I_A$ could in principle be smaller than $A$. However, successive applications of such a partitioning of the index set reduce the size of the set $S_Q$, and thus after not more than $|S_Q|$ steps it produces indecomposable ROSt's. We refer to the partition of the configuration according to this algorithm as the *Q-factorization* of $(\xi, Q)$. The ROSt's which describe the different elements of the above partition are referred here as the *Q-factors* of $(\xi, Q)$. The above considerations readily lead to the following factorization statement.

PROPOSITION 1.7. *Let $(\xi, Q)$ be a ROSt with finite state space that is quasi-stationary under a correlated evolution $\Phi$. Then the Q-factorization yields a time-invariant partition of the point process into indecomposable quasi-stationary ROSt's.*

Ruelle probability cascades provide a class of random overlap structures that are robustly quasi-stationary and indecomposable. These processes, whose structure is motivated by the Parisi hierarchical ansatz [15], were introduced by Ruelle in the context of mean-field spin glass systems, as representing the asymptotics of the models studied by Derrida [9, 16]. An insightful description of RPC through a coalescence process was provided by Bolthausen and Sznitman [8]. In the terminology which is explained in Section 2, a Ruelle probability cascade is a random pair $(\xi, Q)$ where $\xi$ is a Poisson–Dirichlet process $PD(x, 0)$, $x \in (0, 1)$, and $Q$ has a hierarchical structure.

### 1.4. *Statement of the result.* The main result proven here is:

THEOREM 1.8. *Let $\psi \in C^2(\mathbb{R})$ be nonconstant with bounded derivatives, and $\lambda_0 > 0$. If a ROSt is robustly quasi-stationary and ergodic for all multiples $\lambda\psi$ with $|\lambda| < \lambda_0$, then each of its Q-factors is a Ruelle probability cascade. In particular, within each Q-factor, the overlaps $\{q_{ij}\}$ are nonnegative and ultrametric almost surely.*

We note that in the case $\psi(\kappa) = \kappa$ quasi-stationarity for only one $\lambda$ is needed. Moreover, this $\lambda$ can depend on the power of the covariance matrix (cf. Theorem 4.4). It would be interesting to investigate the case of infinite state space, which eludes the approach developed here (see [6] for preliminary results). The proof of the theorem is given in Section 4, through steps which are outlined next.

### 1.5. *Outline of the proof.* The article is organized as follows. Section 2 presents in detail the definition of the Ruelle probability cascades with an emphasis on the quasi-stationarity property. The derivation of the Theorem 1.8 is then divided into four steps.

(i) *The free evolution (Section 3).* As it turns out, it is important to understand the condition of quasi-stationarity of random overlap structures under the "free," that is overlap-independent evolution—for which the Gaussian variables $\kappa$ are taken to be i.i.d. random variables. For point processes without the overlap degrees of freedom it is known that, under a nonlattice condition, quasi-stationarity implies that the law of $\xi$ is a linear superposition of Poisson–Dirichlet processes $PD(x, 0)$ [5, 17]. Random overlap structures include also the overlap matrix. For those we prove:

THEOREM 1.9. *Let $(\xi, Q)$ be a ROSt that is quasi-stationary under the free evolution for a function $\psi$ satisfying Assumption 1.3. There exists a random probability measure $\mu$ on a Hilbert space $\mathcal{H}$ such that conditionally on $\mu$:*

1. *the law of $\xi$ is a linear superposition of $PD(x, 0)$ independent of $Q$;*
2. *$Q$ is directed by $\mu$, that is: for $i \neq j$*

$$q_{ij} \stackrel{\mathcal{D}}{=} (\phi_i, \phi_j)_{\mathcal{H}},$$

*where $(\phi_i, i \in \mathbb{N})$ is i.i.d. $\mu$-distributed.*

The proof is based on the fact that the law of $Q$ must be invariant under the action of any finite permutation (such a random matrix is said to be *weakly exchangeable* [3]). Once this fact is proven, we can apply a variation of de Finetti's theorem due to Dovbysh and Sudakov [10] (see also Hestir [13]). This result states that the law of a weakly exchangeable $\mathbb{N} \times \mathbb{N}$ covariance matrix is directed by a random probability measure $\mu$ similarly to the way that the empirical measure directs the distribution of an exchangeable sequence in the original de Finetti's theorem. The distribution of $\xi$ is obtained through an adaptation of the Ruzmaikina–Aizenman theorem [17] on processes which are quasi-stationary under independent increments [5].

(ii) *Robustness and free evolution (Section 4.1).* Next, we show:

THEOREM 1.10. *If $(\xi, Q)$ is a ROSt that is robustly quasi-stationary for some function $\psi$ satisfying Assumption 1.3, then it is also quasi-stationary under the free evolution.*

In essence, this follows from the fact that, under our assumptions on $Q$, $Q^{*r}$ tends to the identity matrix as $r \to \infty$. By Theorems 1.9 and 1.10 yields a characterization of the law of $\xi$. It remains to single out the distribution of the directing measure. For that we consider the further implications of quasi-stationarity.

(iii) *The directing random overlap structure (Section 4.2).* Pointwise, the directing measure $\mu$, which we recall is defined for each (or almost every) configuration $(\xi, Q)$, evolves under the correlated evolution. However, under the quasi-stationarity assumption its law should be invariant. In Appendix A, we prove that under the assumptions of indecomposability and finiteness of the state space, the directing measure on $\mathcal{H}$ is discrete, of the form $\mu = \sum_{l \in \mathcal{L}} \tilde{\xi}_l \delta_{\phi_l}$ for some countable collection of vectors $(\phi_l, l \in \mathcal{L})$. Hence $\mu$ can be represented by $(\tilde{\xi}, \tilde{Q})$ where $\tilde{\xi}$ are the ordered weights and $\tilde{Q}$ is the Gram matrix of the vectors (with diagonal normalized to 1). $(\tilde{\xi}, \tilde{Q})$ is called the *directing ROSt*. It is rather simple to write down the evolution of this new random overlap structure relying heavily on the fact that $\xi$ is a Poisson–Dirichlet variable. The astounding fact is that $(\tilde{\xi}, \tilde{Q})$ undergoes also an evolution of the form (1.4). However, this evolution is governed by a slightly different function $\psi$ which explicitly depends on $r$. These results are summarized in Proposition 4.3.

(iv) *The induction argument* (*Sections* 4.3 *and* 4.4). The RPC's admit a directing ROSt. Curiously, this structure is again a cascade. This observation plays a role in our analysis, as the proof of Theorem 1.8 proceeds by induction on the cardinality of the state space $S_Q$. When $|S_Q| = 1$, the result of [5, 17] implies that the systems fits the simplest, degenerate, case of the RPC. For the induction step one needs to show that $(\tilde{\xi}, \tilde{Q})$ satisfies the property of indecomposability and robust quasi-stationarity. Indecomposability is evident. Robust quasi-stationarity is trickier as the function of the evolution of $(\tilde{\xi}, \tilde{Q})$ is different for every $r$. This obstacle is circumvented in the linear case by an appropriate scaling between $r$ and the number of time steps of the evolution. The nonlinear case is reduced to the linear one using an argument based on the central limit theorem, where quasi-stationarity is needed for $\lambda$ in a neighborhood of 0. Since $|S_{\tilde{Q}}| = |S_Q| - 1$, the induction proof covers all cases of $|S_Q| < \infty$.

It is conceivable that the robustness assumption of Theorem 1.8 can be weakened. Quasi-stationarity under a single mapping is not enough, as it does not rule out a superposition of independently moving ROSt's which for the given $\psi$ happen to progress at a common velocity. It seems not unreasonable to expect that in situations other than the linear case, $\psi(\kappa) = \kappa$, such degeneracy would be broken by requiring quasi-stationarity for the mapping corresponding to $\psi' \equiv \lambda \psi$ for an open range of values of $\lambda$. It is also not unconceivable that such a condition could directly yield robustness in the sense used in this work. These however are open questions.

## 2. Ruelle probability cascades.

In this section we briefly recount the structure and some of the essential properties of probability cascades which are constructed hierarchically from Poisson–Dirichlet processes. This is done here mainly in order to introduce the notation which is used throughout the paper. The framework, in particular with the insight of Bolthausen and Sznitman [8], allows an extension of the definition to cascades with an infinite number of levels of splitting. In this article, we focus on the finite case. In discussing the quasi-stationarity we formulate a slight extension of the known stationarity property of the Poisson–Dirichlet processes.

The probabilistic cascade models were introduced by Ruelle [16] as a way of giving expression to Derrida's asymptotic calculations of the $N \to \infty$ limits of a class of finite models [9]. Often they are referred to as GREM and REM, for (*generalized-*) *random energy models*, though this terminology runs the risk of confusing the finite models with their continuum reformulation.

A Poisson point process $\eta$ on a Polish space $S$ with intensity measure $\nu$ is the integer-valued random measure whose finite-dimensional distributions are of the form

$$\mathbb{P}(\eta(A_1) = m_1, \ldots, \eta(A_n) = m_n) = \prod_{k=1}^{n} \frac{\nu(A_k)^{m_k}}{m_k!} e^{-\nu(A_k)}$$

for disjoint measurable sets $A_k$ and $m_k \in \mathbb{N}$. We are particularly interested in the Poisson processes on $(0, \infty)$ with intensity measure $\nu(ds) = x s^{-x-1}\, ds$ for some $x \in (0, 1)$, whose tail is given by the power law $\nu([s, \infty)) = s^{-x}$. (Some relevant basic properties of this process are mentioned in [2].)

DEFINITION 2.1. A Poisson–Dirichlet variable $PD(x, 0)$, $x \in (0, 1)$, on $(P_m, \mathcal{F}_m)$ is the random mass-partition defined as

$$\xi := \left( \frac{\eta_i}{\sum_j \eta_j}, i \in \mathbb{N} \right),$$

where $(\eta_i, i \in \mathbb{N})$ are the ordered positions of the atoms of a Poisson process on $(0, \infty)$ with intensity measure $x s^{-x-1}\, ds$.

It should be noted that the normalization is by a finite factor, since $\sum_i \eta_i < \infty$ a.s. for $x \in (0, 1)$. Furthermore, this factor can be almost surely determined from the normalized weights of the mass partition, as in the following explicit statement whose proof is elementary (see, e.g., [2, 7]).

PROPOSITION 2.2. For $\eta = (\eta_n, n \in \mathbb{N})$ a Poisson process, as above, and $\xi = (\xi_n, n \in \mathbb{N})$ the corresponding $PD(x, 0)$ variable:

$$\lim_{n \to \infty} n^{1/x} \xi_n \overset{a.s.}{=} \frac{1}{\zeta},$$

where $\zeta = \sum_i \eta_i$ is the normalizing factor relating $\xi$ and $\eta$.

In particular: (i) $\zeta$ is measurable with respect to $\xi$, (ii) the process $(\zeta \xi_n, n \in \mathbb{N})$ is Poisson with intensity measure $x s^{-x-1}\, ds$, and (iii) the laws of Poisson–Dirichlet variables with distinct parameters are mutually singular.

2.1. *Definition of the hierarchical structure.* The family of Ruelle probability cascades or RPC's, including the ones with continuous branching, is parametrized by the space of probability measures on $[0, 1]$. From this perspective, the cascades with finite number of splittings are in correspondence with the discrete probability measures with finite number of atoms, one of which is located at 1.

More precisely, for each *number of level splittings*, $k \in \mathbb{N}$, the distribution of the cascade is determined by the $2k$ parameters:

- the values taken by the overlaps

$$0 \leq q_1 < \cdots < q_k < q_{k+1} = 1,$$

- the parameters of the Poisson processes

$$0 < x_1 < \cdots < x_k < x_{k+1} = 1.$$

These $2k$ parameters define a piecewise constant distribution function

$$(2.1) \qquad x(q) := \begin{cases} 0, & \text{for } q \in [0, q_1), \\ x_l, & \text{for } q \in [q_l, q_{l+1}). \end{cases}$$

This corresponds to a discrete probability measure on $[0,1]$ with $k+1$ atoms.

The cascade parametrized by $x(q)$ is constructed as follows. Let $\alpha := (\alpha_1, \ldots, \alpha_k) \in \mathbb{N}^k$. It is convenient to define $\alpha(l)$ as the truncation of $\alpha$ up to the $l$th component, that is, $\alpha(l) := (\alpha_1, \ldots, \alpha_l)$. By convention, $\alpha(0) := 0$. The collection of $\alpha \in \mathbb{N}^k$ and their truncations is naturally represented as a rooted tree with root 0, vertices $\alpha(l)$ and leaves $\alpha$. Edges are drawn between $\alpha(l+1)$ and its truncation $\alpha(l)$, for every $\alpha \in \mathbb{N}^k$ and $0 \leq l \leq k$. The Poisson processes constituting the cascade are indexed by the vertices of this tree. Precisely, to each vertex $\alpha(l)$ and for every $0 \leq l \leq k-1$, we associate an independent Poisson process $\eta^{\alpha(l)}$ with intensity measures $x_{l+1} s^{-x_{l+1}-1} ds$. The law of the RPC, seen as a ROSt, expresses this hierarchy of processes.

We first define the point process $\eta$ with atoms indexed by $\alpha \in \mathbb{N}^k$

$$(2.2) \qquad \eta = (\eta_\alpha, \alpha \in \mathbb{N}^k) = (\eta^0_{\alpha_1} \eta^{\alpha(1)}_{\alpha_2} \cdots \eta^{\alpha(k-1)}_{\alpha_k}, \alpha \in \mathbb{N}^k),$$

where $\eta^{\alpha(l)}_{\alpha_{l+1}}$ represents the position of the $\alpha_{l+1}$th atom of the Poisson process $\eta^{\alpha(l)}$. The point process $\eta = (\eta_\alpha, \alpha \in \mathbb{N}^k)$ retains some properties of the Poisson processes constituting it (see, e.g., [2]). In particular it can be ordered as it is summable. Hence we construct the random mass-partition

$$(2.3) \qquad \xi := \left( \frac{\eta_i}{\sum_j \eta_j}, j \in \mathbb{N} \right).$$

As the atoms of $\xi$ are ordered, there exists a random bijection $\Pi : \mathbb{N} \to \mathbb{N}^k$ such that $\Pi(i) = \alpha$ if the $i$th point of $\xi$ has address $\alpha$. This map induces a hierarchy of equivalence relations: for each $l \in \{0, \ldots, k\}$

$$i \sim_l j \quad \Longleftrightarrow \quad [\Pi(i)](l) = [\Pi(j)](l).$$

In other words, $i$ and $j$ are equivalent under $\sim_l$ if and only if $\xi_i$ and $\xi_j$ share a common ancestor at level $l$ in the tree indexed by $\alpha$. We write $\Gamma_l$ for the partition of $\mathbb{N}$ induced by this equivalence relation. It is clear from the definition that $i \sim_{l'} j$ for every $l' < l$ whenever $i \sim_l j$, or equivalently $\Gamma_{l'} \subset \Gamma_l$. Let $l(i,j)$ be the greatest $l$ for which $i$ and $j$ belong to the same block of $\Gamma_l$. The *overlap* between the $i$th point and $j$th point of $\xi$ is defined as

$$(2.4) \qquad q_{ij} := q_{l(i,j)+1},$$

where $0 \leq q_1 < \cdots < q_k < q_{k+1} = 1$ are the overlap parameters of the cascade. The random matrix $Q$ hereby constructed is clearly real symmetric with $q_{ii} = 1$. It is also positive definite. Indeed, define the Gaussian field indexed by $\mathbb{N}^k$

$$(2.5) \qquad \kappa_\alpha := \sum_{l=0}^{k} \Delta\kappa_{\alpha(l)},$$

where $\Delta\kappa_{\alpha(l)}$ are independent centered Gaussians with variance $q_{l+1} - q_l$ for all $\alpha$ and $0 \leq l \leq k$. It is easily checked that the field $(\kappa_i, i \in \mathbb{N})$ has covariance matrix $Q$ thereby proving the positivity.

DEFINITION 2.3 (*Ruelle probability cascade*). A $k$-level Ruelle probability cascade or RPC with parameter $x(q)$ (or equivalently parameters $0 < x_1 < \cdots < x_k < x_{k+1} = 1$ and $0 \leq q_1 < \cdots < q_k < q_{k+1} = 1$) is the ROSt $(\xi, Q)$ defined by (2.2) and (2.4).

Recall that the $r$th power $Q^{*r}$ (in the sense of the entrywise product) is positive definite whenever $Q$ is. In fact, if $(\xi, Q)$ is a RPC then so is $(\xi, Q^{*r})$. An interesting property coming from the hierarchical organization of the cascade is the fact that the latter holds for any monotone increasing function of the entries.

PROPOSITION 2.4. *Let* $F : [0, 1] \to [0, 1]$ *be a strictly increasing function such that* $F(1) = 1$. *If* $(\xi, Q)$ *is a* $k$-level RPC, *then so is* $(\xi, F(Q))$ *where* $F(Q) := \{F(q_{ij})\}$.

PROOF. It suffices to take $\Delta\kappa_{\alpha(l)}$ i.i.d. centered Gaussian with variance $F(q_{l+1}) - F(q_l)$ in (2.5). □

Bolthausen and Sznitman (see Theorem 2.2 in [8]) introduced a simple description of the law of the cascade which avoids the reordering of the index $\alpha$. The following formulation of their result will be used in the proof of Theorem 1.8. Note that the 1-level RPC with parameter $x_1$ and $q_1$ is simply a Poisson–Dirichlet variable $PD(x_1, 0)$ together with the deterministic overlap matrix with nondiagonal entries $q_1$.

THEOREM 2.5 [8]. *Let* $(\xi, Q)$ *be a* $k$-level RPC, $k > 1$, *with parameters* $x_1, \ldots, x_k$ *and* $q_1, \ldots, q_k$. *Then:*

1. $\xi$ and $Q$ are distributed independently, and

2. $\xi$ is a $PD(x_k, 0)$ variable.

3. The distribution of $Q$ can be constructed as follows. Let $(\tilde{\xi}, \widetilde{Q})$ be a $(k-1)$-level probability cascade with parameters $\frac{x_1}{x_k}, \ldots, \frac{x_{k-1}}{x_k}$ and $\frac{q_1}{q_k}, \ldots, \frac{q_{k-1}}{q_k}$. Conditionally on $(\tilde{\xi}, \widetilde{Q})$, let $(i^*, i \in \mathbb{N})$ be i.i.d. $\tilde{\xi}$-distributed random integers, then

$$(2.6) \qquad q_{ij} \stackrel{\mathcal{D}}{=} q_k \tilde{q}_{i^* j^*},$$

where $\stackrel{\mathcal{D}}{=}$ means equality of distributions.

We stress the fact, used in the above, that any proper mass-partition $\xi = (\xi_i, i \in \mathbb{N})$ can be seen as a probability measure on $\mathbb{N}$ where the integer $i$ is sampled with probability $\xi_i$. The auxiliary structure $(\tilde{\xi}, \widetilde{Q})$ appearing in the construction is said to be the *directing probability cascade*.

2.2. *Quasi-stationarity of the cascade.* As was stated previously, the Ruelle probability cascades are examples of random overlap structures that are indecomposable and robustly quasi-stationary under the correlated evolution, see, for example, [2, 4]. In discussing this claim, it is useful to bear in mind the following result on the shift of a marked Poisson–Dirichlet variable.

DEFINITION 2.6. Let $\mathcal{C}$ be a Polish space equipped with the probability measure $\mu$, $\eta$ a Poisson process on $\mathbb{R}^+$, and $\xi$ be Poisson–Dirichlet process $PD(x, 0)$. The following are referred to as the corresponding $(\mathcal{C}, \mu)$-marked processes:

$$\eta^{\mathcal{C}} := ((\eta_i, c_i), i \in \mathbb{N}), \qquad \xi^{\mathcal{C}} := ((\xi_i, c_i), i \in \mathbb{N}).$$

The variable $c_i$ is called the mark of the point $i$.

The following statement provides a slight generalization of Proposition 3.1 in [17].

LEMMA 2.7. *Let $\eta$ be a Poisson process on $(0, \infty)$ with intensity measure $x s^{-x-1} ds$ and $W_\cdot(\cdot) \colon \mathcal{C} \times \mathbb{R} \to (0, \infty)$ be a jointly measurable mapping such that*

$$\mathcal{N} := \int_{\mathcal{C}} \int_{\mathbb{R}} W_c^x(\kappa) \mu(dc) \rho_c(d\kappa) < \infty$$

*for some probability measure $\mu(dc)\rho_c(d\kappa) \equiv \rho(dc\,d\kappa)$ on the product space $(\mathcal{C} \times \mathbb{R})$. Consider a $(\mathcal{C} \times \mathbb{R}, \rho)$ marking of $\eta$. Then the shifted process*

$$(\eta_i W_{c_i}(\kappa_i), c_i)$$

*is a $(\mathcal{C}, \tilde{\mu})$-marked Poisson process with intensity measure $\mathcal{N} x s^{-x-1} ds$ and*

$$\tilde{\mu}(dc) := \frac{\int_{\mathbb{R}} W_c^x(\kappa) \rho_c(d\kappa)}{\mathcal{N}}.$$

The statement has an elementary proof using the characteristic functional of Poisson processes, or alternatively it also follows by a direct adaptation of the proof of Proposition 3.1 in [17]. It has the following immediate consequence.

COROLLARY 2.8.  *In the notation of Lemma 2.7, consider a $(\mathcal{C} \times \mathbb{R}, \rho)$-marking of a $PD(x,0)$ variable. Then, the shifted process*

$$\left( \frac{\xi_i W_{c_i}(\kappa_i)}{\sum_j \xi_j W_{c_j}(\kappa_j)}, c_i \right)$$

*is a $(\mathcal{C}, \tilde{\mu})$-marked $PD(x,0)$.*

We will usually deal with the uncoupled case where $\rho_c(d\kappa)$ does not depend on the mark $c$. Note that if moreover $W_c \equiv W$, then $\tilde{\mu} = \mu$.

For the latter case, let us describe the distribution of the past increments. Let $\pi$ be the random permutation of $\mathbb{N}$ induced by the shift, that is: $\pi(i) = j$ if and only if $\frac{\xi_i W(\kappa_i)}{\sum_k \xi_k W(\kappa_k)}$ is the $j$th point of the evolved process. We consider $\kappa'_i := \kappa_{\pi^{-1}(i)}$. This variable represents the past increment of the point that made it to the $i$th position. The distribution of $(\kappa'_i, i \in \mathbb{N})$ is computed by simply considering $\kappa_i$ as a mark. Lemma 2.7 shows that the past increments are i.i.d. with distribution

$$\frac{W(\kappa)^x \rho(d\kappa)}{\mathcal{N}}$$

and independent of $(\frac{\xi_i W(\kappa_i)}{\sum_j \xi_j W(\kappa_j)}, i \in \mathbb{N})$. We note for further references that the above remark is easily extended to the case where the points of the process are incremented by $T$ independent steps. In this case, the past steps are again independent.

THEOREM 2.9.  *Let $(\xi, Q)$ be a $k$-level RPC for some $k \in \mathbb{N}$. Then $(\xi, Q)$ is robustly quasi-stationary and indecomposable under the correlated evolution (1.4) for any function $\psi$ satisfying Assumption 1.3.*

PROOF.  It suffices to prove quasi-stationarity for $r = 1$. Indeed, $(\xi, Q^{*r})$ is also a cascade by Proposition 2.4 hence it is quasi-stationary for the Gaussian field with covariance $Q^{*r}$. As $Q = F(Q^{*r})$ for the increasing function $F(q) = q^{1/r}$, we conclude that $(\xi, Q)$ is quasi-stationary for the field with covariance $Q^{*r}$.

The proof is by induction on $k$. In the case $k = 1$, the matrix $Q$ is simply $q_{ij} = \delta_{ij} + (1 - \delta_{ij})q$ for some $0 \le q < 1$. In particular, it is invariant under the correlated evolution. The field $\kappa$ can then be represented as

$$\kappa_i = \kappa^c + \kappa_i^f,$$

where $\kappa^c$ is a Gaussian with variance $q$ independent of the i.i.d. Gaussians $\kappa_i^f$ of variance $1 - q$. Conditioning on $\kappa^c$, we have that the law of $\xi$ is the same independently of $\kappa^c$ under the evolution by Corollary 2.8 taking $W(\kappa_i^f) = e^{\psi(\kappa^c + \kappa_i^f)}$. The case $k = 1$ is proven.

We now assume that all $(k-1)$-level cascades are quasi-stationary. By Proposition 2.5, the distribution of a $k$-level cascade is such that $\xi$ is $PD(x_k, 0)$ independent of $Q$ and that there exists a $(k-1)$-level cascade $(\tilde{\xi}, \tilde{Q})$ such that

$$(2.7) \qquad\qquad q_{ij} \overset{\mathcal{D}}{=} q_k \tilde{q}_{i^* j^*},$$

where $(i^*, i \in \mathbb{N})$ are independent $\tilde{\xi}$-distributed random integers. In particular, we can write

$$\kappa_i = \kappa_{i^*}^c + \kappa_i^f,$$

where the $\kappa_i^f$ are i.i.d. Gaussians with variance $1 - q_k$ and $(\kappa_l^c, l \in \mathbb{N})$ is a Gaussian field with variance $q_k \tilde{Q}$.

Let us condition on $(\tilde{\xi}, \tilde{Q})$. Clearly, the process $(\xi, Q)$ is equivalent to the $(\mathbb{N}, \tilde{\xi})$-marking of a $PD(x_k, 0)$ variable

$$\xi^{\mathbb{N}} = ((\xi_i, i^*), i \in \mathbb{N})$$

with the mark probability $\tilde{\xi} = (\tilde{\xi}_l, l \in \mathbb{N})$ by (2.7). Quasi-stationarity will be proven if we show that the evolved process is still a $\mathcal{C}$-marked $PD(x_k, 0)$ with the same mark probability (in law). Conditionally on $(\kappa_l^c, l \in \mathbb{N})$, the correlated evolution of (1.4) corresponds to the evolution of the marking in Corollary 2.8 with function $W_l(\kappa_i^f) = e^{\psi(\kappa_l^c + \kappa_i^f)}$ for $l \in \mathbb{N}$. Therefore, the evolved process is still a marked $PD(x_k, 0)$. It remains to show that the law of the evolved directing cascade is the same as $(\tilde{\xi}, \tilde{Q})$. The mark probability after evolution is by Corollary 2.8

$$(2.8) \qquad\qquad \left( \frac{\tilde{\xi}_l e^{\psi_{x_k, q_k}(\kappa_l^c)}}{\sum_{l'} \tilde{\xi}_{l'} e^{\psi_{x_k, q_k}(\kappa_{l'}^c)}} \right),$$

where $e^{\psi_{x, \rho}(y)} := \int_{\mathbb{R}} \frac{e^{-z^2/2}}{\sqrt{2\pi}} e^{x\psi(y + z\sqrt{1-\rho})} \, dz$. Note that this function is well defined by Assumption 1.3. The reordering of the evolved mark probability in (2.8) induces a permutation of $\tilde{Q}$

$$\tilde{Q} \mapsto \tilde{\pi} \tilde{Q} \tilde{\pi}^{-1},$$

where $\tilde{\pi}(l) = k$ if and only if $\tilde{\xi}_l$ becomes the $k$th highest mark weight in (2.8). Hence the evolution of $(\tilde{\xi}, \tilde{Q})$ induced by the evolution of $(\xi, Q)$ is exactly a correlated evolution with function $\psi_{x_k, q_k}$. By induction, the law of $(\tilde{\xi}, \tilde{Q})$ is invariant under this correlated evolution. Quasi-stationarity of $(\xi, Q)$ is proven.

Indecomposability is also proven by induction on $k$. The case $k = 1$ is obvious. Assume every $(k-1)$-level cascade is indecomposable. So is any $k$-level cascade as

$$S_Q(i) = q_k S_{\tilde{Q}}(i^*) = q_k S_{\tilde{Q}}$$

from (2.7) and the induction hypothesis.  □

The existence of the directing cascade $(\tilde{\xi}, \tilde{Q})$ and its evolution as a marking of $\xi$ are two important elements that will reappear in the characterization of quasi-stationary laws. We are now ready to turn to the implications of quasi-stationarity of a ROSt.

**3. Quasi-stationarity under the free evolution.**   In this section, we study the overlap-independent or free evolution of random overlap structures. The evolution is by the mapping (1.4) where the variables $(\kappa_i, i \in \mathbb{N})$ are simply i.i.d. standard Gaussians, whose covariance matrix is the identity. We denote this mapping by $\Phi_\infty$, and will write $\mathbb{P}_\infty$ for the joint law of the ROSt and the Gaussian field. The goal is to characterize the law of random overlap structures that are quasi-stationary under the free evolution for a given function $\psi$. The full description is stated as Theorem 1.9 in the Introduction.

3.1. *Characterization theorem.*   To achieve the desired characterization, it is necessary to assume that the evolution is nonlattice. The next theorem was proven in [5] based on the argument of Ruzmaikina and Aizenman [17], whose result is stated under somewhat more stringent assumptions.

THEOREM 3.1.   *Let $\xi$ be a random mass-partition such that $\xi_1 \neq 1$ a.s. Let $(W_i, i \in \mathbb{N})$ be an i.i.d. sequence of positive numbers whose distribution has a density on $(0, \infty)$ and for which $\mathbb{E}[W^\lambda] < \infty$ for any $\lambda \in \mathbb{R}$. If the law of $\xi$ is invariant under the evolution*

$$(3.1) \qquad (\xi_i, i \in \mathbb{N}) \mapsto \left( \frac{\xi_i W_i}{\sum_j \xi_j W_j}, i \in \mathbb{N} \right)_\downarrow,$$

*then it is a linear superposition of Poisson–Dirichlet distributions $PD(x, 0)$, $x \in (0, 1)$.*

The free evolution corresponds to the choice $W_i = e^{\psi(\kappa_i)}$. Clearly, the assumptions of Theorem 3.1 are fulfilled when $\psi$ satisfies Assumptions 1.3. We conclude that the marginal distribution of $\xi$ must be a superposition of Poisson–Dirichlet variables.

The key point of the proof of Theorem 1.9 is the fact that $Q$ must be weakly exchangeable given $\xi$ (i.e., its law must be invariant under any finite permutation of the indices). The properties of weakly exchangeable overlap matrices needed for the section are detailed in Appendix A. The proof of the weak exchangeability of the overlap matrix is given in the next section (Proposition 3.3). We first complete the proof of the theorem assuming this fact.

PROOF OF THEOREM 1.9.    It will be proven in Proposition 3.3 that $Q$ is weakly exchangeable given $\xi$. The characterization of weakly exchangeable covariance matrices of Dovbysh and Sudakov (Theorem A.1) guarantees the existence of a random probability measure $\mu$ on $\mathcal{H}$ such that the law of $Q$ has the desired form conditionally on $\mu$. Moreover, given $\mu$, the law of $Q$ does not depend on $\xi$ thereby proving the independence between $\xi$ and $Q$.

It remains to prove that the law of $\xi$ conditionally of $\mu$ is invariant under the free evolution. The fact that it must be a superposition of $PD(x, 0)$ will then follow by Theorem 3.1. We write $\mu_Q$ for the probability measure making the dependence on $Q$ explicit. We denote the overlap matrix of $\Phi_\infty(\xi, Q)$ by $Q'$. The claim will be proven if we show that the directing measure of $Q$ is a $\Phi_\infty$-invariant, that is: for all measurable $A \subseteq \mathcal{H}$

(3.2) $$\mathbb{E}_\infty[\mu_{Q'}(A)|\mu_Q] = \mu_Q(A).$$

The probability $\mu_{Q'}(A)$ is the probability that the first point (or any point) of $\Phi_\infty(\xi, Q)$ has its mark $\phi_1$ in $A$

$$\mathbb{E}_\infty[\mu_{Q'}(A)|\mu_Q] = \mathbb{P}_\infty(\phi_{\pi^{-1}(1)} \in A|\mu_Q),$$

where $\pi$ is the reshuffling induced by the evolution. By summing over all possibilities of leading points, we get

$$\mathbb{E}_\infty[\mu_{Q'}(A)|\mu_Q] = \sum_j \mathbb{P}_\infty(\phi_j \in A, \pi(j) = 1|\mu_Q).$$

As the increments $\kappa$ (and thus the reshuffling) are independent of the marks $\phi_j$, the r.h.s becomes

$$\sum_j \mathbb{P}(\phi_j \in A|\mu_Q)\mathbb{P}_\infty(\pi(j) = 1|\mu_Q) = \mu_Q(A),$$

where we have used the fact that $\mathbb{P}(\phi_j \in A|\mu_Q) = \mu_Q(A)$ for all $j$. The claim is proven.    □

We remark that the theorem is sharp in the sense that any random overlap structure where $\xi$ is a superposition of $PD(x,0)$ and $Q$ is constructed from a probability measure on $\mathcal{H}$ is quasi-stationary under the free evolution. This is a direct consequence of Corollary 2.8 for the choice $\mathcal{C} = \mathcal{H}$ when the shift function does not depend on the mark. In other words, quasi-stationarity under the free evolution does not constrain the form of the directing measure $\mu$. However, we will see that indecomposability and the finiteness of the state space restrict the possible $\mu$ in a dramatic way.

3.2. *Weak exchangeability of the overlap matrix.* We shall now prove that weak exchangeability of the overlap matrix is a necessary condition for quasi-stationarity under the free evolution. The proof is obtained by looking at the distribution of the past increments given the present positions of the points.

Consider independent copies $\kappa'(t)$, $0 \leq t \leq T-1$, of the field $\kappa$ composed of i.i.d. standard Gaussians. Let $(\xi', Q')$ be a ROSt and define the evolved process $(\xi, Q)$ after the $T$-step incrementation

$$(\xi, Q) := \Phi_{\psi(\kappa'(T-1))} \circ \cdots \circ \Phi_{\psi(\kappa'(0))}(\xi', Q').$$

We define $\pi_T$ as the reshuffling after these $T$ steps, that is, $\pi_T(i) = j$ if and only if

$$\xi_j = \frac{\xi'_i e^{\sum_{t=0}^{T-1} \psi(\kappa'_i(t))}}{\sum_k \xi'_k e^{\sum_{t=0}^{T-1} \psi(\kappa'_k(t))}}.$$

The past Gaussian field at time $-t$, $-T \leq -t \leq -1$, of the $i$th point of $(\xi, Q)$ is then

$$(3.3) \qquad \kappa_i(-t) := \kappa'_{\pi_T^{-1}(i)}(T-t).$$

LEMMA 3.2. *Let $(\xi, Q)$ be a ROSt that is quasi-stationary under the free evolution for some function $\psi$ satisfying Assumption 1.3. Consider the past fields as defined in (3.3) for some $T \in \mathbb{N}$. There exists a $\mathcal{F}_{os}$-measurable parameter $x$ such that, conditionally on $x$, the following hold:*

1. *The fields $\kappa(t)$, $-T \leq -t \leq -1$, are independent of each other and of $(\xi, Q)$.*
2. *For fixed $t$, the variables $(\kappa_i(-t), i \in \mathbb{N})$ are i.i.d. with distribution*

$$\frac{e^{x\psi(z)} \nu(dz)}{\int_{\mathbb{R}} e^{x\psi(z')} \nu(dz')},$$

*where $\nu$ is the standard Gaussian measure.*

PROOF.   We know from Theorem 1.9 that the law of $\xi$ is a superposition of $PD(x,0)$. As the increments are independent of the overlaps, we can apply directly Corollary 2.8 and the remark following it to get the distribution of the fields $\kappa(-t)$, $-T \le -t \le -1$.   □

Note that the field $(\kappa(-t), t \in \mathbb{N})$ is well defined as i.i.d. copies of $\kappa(-1)$. In the rest of this section, we will write $\mathbb{P}_\infty$ for the distribution of $(\xi, Q)$ together with the law of the past and future Gaussian fields $(\kappa(-t), t \in \mathbb{Z})$.

PROPOSITION 3.3.   *Let $(\xi, Q)$ be a ROSt that is quasi-stationary under free evolution for some function $\psi$ satisfying Assumption 1.3. Then $Q$ is weakly exchangeable conditionally on $\xi$, that is: for any permutation $\tau$ of a finite number of elements of $\mathbb{N}$*

$$\tau Q \tau^{-1} \stackrel{\mathcal{D}}{=} Q.$$

PROOF.   Let us fix $\tau$, a permutation of a finite number of elements of $\mathbb{N}$. Let $n$ be the maximum over $\mathbb{N}$ of the elements for which $\tau$ is not the identity. We write $Q_n$ for the $n \times n$ matrix of the overlaps of the first $n$ points. Let $I_i$, $i \in \mathbb{N}$, be a collection of intervals of $[0,1]$. We will write $\{\xi \in I\}$ for the event $\{\xi_i \in I_i, i \in \mathbb{N}\}$. We need to prove that

(3.4)          $$\mathbb{P}(Q_n = \tau R_n \tau^{-1}, \xi \in I) = \mathbb{P}(Q_n = R_n, \xi \in I)$$

for any realization $R_n$ of $Q_n$.

Let $\pi_{n,-T} \colon \{1, \ldots, n\} \to \{1, \ldots, n\}$ be the reshuffling of the first $n$ points from time $-T$ to time $0$, that is: $\pi_{n,-T}(i) = j$ if the $j$th point at time $0$ was the $i$th point at time $-T$ among these $n$ points (see Figure 1). We write

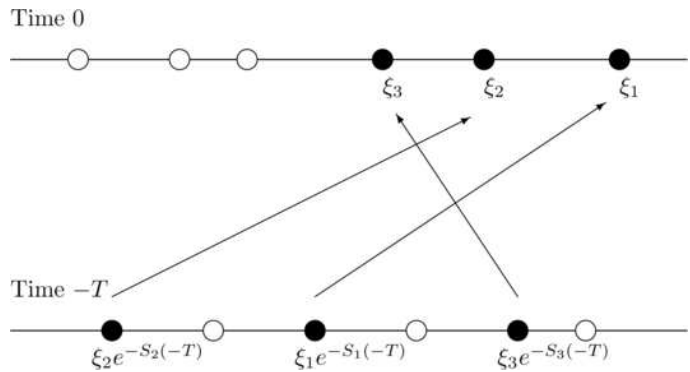

FIG. 1.   *An illustration of the permutations which occur through the time evolution. For each $n \in \mathbb{N}$, the permutation $\pi_{n,-T}$ describes the effect of the evolution from time $-T$ on the $n$ leading points of time $0$. For example, in the depicted case $\pi_{3,-T}(2) = 1$.*

$Q_{n,-T}$ for the overlap matrix of the ordered points at time $-T$ that become the first $n$ at time $0$

$$Q_{n,-T} = \pi_{n,-T}^{-1} Q_n \pi_{n,-T}.$$

By decomposing $\mathbb{P}(Q_n = \tau R_n \tau^{-1}, \xi \in I)$ over all possible $\pi_{n,-T}$ and using conditioning, we get

$$\sum_{\rho' \in S_n} \mathbb{P}_\infty(\pi_{n,-T} = \rho' \tau^{-1} | Q_{n,-T} = \rho' R_n \rho'^{-1}, \xi \in I)$$

$$\times \, \mathbb{P}_\infty(Q_{n,-T} = \rho' R_n \rho'^{-1}, \xi \in I).$$

Equation (3.4) will follow if, for any permutation $\rho \in S_n$ and any realization $R_n'$ of $Q_n$,

$$(3.5) \qquad \mathbb{P}_\infty(\pi_{n,-T} = \rho | Q_{n,-T} = R_n', \xi \in I) \to \frac{1}{n!}$$

as $T \to \infty$.

From the definition of $\pi_{n,-T}$, the l.h.s. of the above equals

$$\mathbb{P}_\infty(\xi_{\rho(1)} e^{-S_{\rho(1)}(-T)} \geq \cdots \geq \xi_{\rho(n)} e^{-S_{\rho(n)}(-T)} | Q_{n,-T} = R_n', \xi \in I),$$

where $S_j(-T) := \sum_{t=1}^T \kappa_j(-t)$ (see Figure 1); and by taking the logarithm

$$\mathbb{P}_\infty\left(\Delta_{\rho(i)}(T) \geq \log \frac{\xi_{\rho(i+1)}}{\xi_{\rho(i)}}, i = 1, \dots, n | Q_{n,-T} = R_n', \xi \in I\right),$$

where we wrote $\Delta_i(T)$ for $S_{i+1}(-T) - S_i(-T)$. By Lemma 3.2, the $n$-dimensional variable $(\Delta_1(T), \dots, \Delta_n(T))$ is a sum of $T$ i.i.d. $n$-dimensional vectors conditionally on the $\mathcal{F}_{os}$-measurable parameter $x$. It is also independent of $\xi$ and the overlaps. We know from the specific form of its distribution that $\mathbb{E}[\Delta_i] = 0$ for all $i$ and $\mathbb{E}[\Delta_i \Delta_j] = CT\delta_{ij}$ for some finite $C$. Dividing all inequalities in the above by $\sqrt{CT}$ does not change the probability and we get

$$\mathbb{P}_\infty(\pi_{n,-T} = \rho | Q_{n,-T} = R_n', \xi \in I)$$

$$= \int_{(0,1)} \lambda(dx) \mathbb{P}_\infty\left(\frac{1}{\sqrt{CT}} \Delta_{\rho(i)}(T)\right.$$

$$\left. \geq \frac{1}{\sqrt{CT}} \log \frac{\xi_{\rho(i+1)}}{\xi_{\rho(i)}}, i = 1, \dots, n | x, \xi \in I\right)$$

for some probability measure $\lambda$ on $(0,1)$. We can take the limit $T \to \infty$ inside the integral by dominated convergence and use the central limit theorem. As the positions $\xi$ are fixed, the r.h.s. of the inequalities goes to 0 and the probability converges to $1/n!$ irrespective of $x$. Equation (3.5) is proven and the proposition follows. $\square$

**4. Robust quasi-stationarity and hierarchical structure.** We now turn to random overlap structures that are quasi-stationary under the correlated evolutions $\Phi_r$ for an infinite number of $r \in \mathbb{N}$. The complete description is stated in Theorem 1.8 for the case where the invariance holds for multiples of a smooth function $\psi$, and the overlap matrix has finite state space.

4.1. *Robustness and free evolution.* We first show that robust quasi-stationarity under correlated evolution implies quasi-stationarity under the free evolution $\Phi_\infty$ (Theorem 1.10). Using results of the previous section, this yields a precise form for the robustly quasi-stationary laws. The proof relies on the next lemma. It asserts that with large probability the particles that are located in $[\delta, 1]$ after evolution came from the first $N$ before evolution. We state the result for $T$ independent step of the evolution, as it will be needed later in the proof. We will write $\kappa(t)$, $0 \le t \le T-1$, for $T$ independent copies of the Gaussian field.

LEMMA 4.1. *Let $(\xi, Q)$ be a ROSt such that $\xi$ is a.s. proper, and $\Phi_r$ a correlated evolution for some $r \in \mathbb{N}$ and a function $\psi$ satisfying Assumption 1.3. Let $A_{N,T,\lambda,\delta}$ be the event that a particle initially beyond the $N$th after $T$ steps of the evolution ends within the interval $[\delta, 1]$:*

$$A_{N,T,\lambda,\delta} = \bigcup_{i>N} \left\{ \frac{\xi_i e^{\sum_{t=0}^{T-1} \lambda\psi(\kappa_i(t))}}{\sum_j \xi_j e^{\sum_{t=0}^{T-1} \lambda\psi(\kappa_j(t))}} > \delta \right\}.$$

*Then*

$$\mathbb{P}_r(A_{N,T,\lambda,\delta}) \le \frac{1}{\delta}(g(2\lambda)g(-2\lambda))^{T/2} \sum_{i>N} \mathbb{E}[\xi_i]$$

*for $g(\lambda) = \int_{\mathbb{R}} \frac{e^{-z^2/2}}{\sqrt{2\pi}} e^{\lambda\psi(z)}\, dz$.*

PROOF. Fix $\delta > 0$. Markov's inequality yields the elementary bound

$$\mathbb{P}_r\left( \bigcup_{i>N} \left\{ \frac{\xi_i e^{\sum_{t=0}^{T-1} \lambda\psi(\kappa_i(t))}}{\sum_j \xi_j e^{\sum_{t=0}^{T-1} \lambda\psi(\kappa_j(t))}} > \delta \right\} \right) \le \frac{1}{\delta}\mathbb{E}_r\left[ \frac{\sum_{i>N} \xi_i e^{\sum_{t=0}^{T-1} \lambda\psi(\kappa_i(t))}}{\sum_j \xi_j e^{\sum_{t=0}^{T-1} \lambda\psi(\kappa_j(t))}} \right].$$

We use Cauchy–Schwarz inequality followed by two applications of Jensen's inequality with the functions $f(y) = y^2$ and $f(y) = 1/y^2$ on the r.h.s. to get the upper bound

$$\frac{1}{\delta}(g(2\lambda)g(-2\lambda))^{T/2} \sum_{i>N} \mathbb{E}[\xi_i].$$

Note that the condition $\sum_i \xi_i = 1$ a.s. is needed to apply Jensen's inequality. □

The important feature of the result is that the estimate is uniform in $r$ for all evolutions $\Phi_r$. This allows us to take a rigorous limit $r \to \infty$ under which $Q^{*r}$ tends to the identity matrix thereby decorrelating the evolution.

PROOF OF THEOREM 1.10.   As $(\xi, Q)$ is quasi-stationary, we can assume that $\sum_i \xi_i = 1$ and that there are an infinite number of atoms with $\xi_i > 0$. Recall from Proposition 1.2 that the law of a ROSt is determined by the class of continuous functions that depend on a finite number of points. Let $f : \Omega_{os} \to \mathbb{R}$ be such a function depending on the first $n$ points for some $n \in \mathbb{N}$. We will write $f(\xi, Q) = f(\xi_1, \ldots, \xi_n; Q_n)$ for the explicit dependence where $Q_n = \{q_{ij}\}_{1 \leq i,j \leq n}$. Robust quasi-stationarity implies that for all $r \in \mathbb{N}$

$$(4.1) \qquad \mathbb{E}_r[f(\Phi_r(\xi, Q))] = \mathbb{E}[f(\xi, Q)].$$

The theorem will be proven by showing that the limit $r \to \infty$ of the above equation yields the equality for the free evolution, that is

$$(4.2) \qquad \mathbb{E}_\infty[f(\Phi_\infty(\xi, Q))] = \mathbb{E}[f(\xi, Q)].$$

For $\delta > 0$, we define the function $f$ with a cutoff

$$f_\delta(\xi_1, \ldots, \xi_n; Q_n) := f(\xi_1, \ldots, \xi_n; Q_n) \chi_{\{\xi_n \geq \delta\}},$$

where $\chi_A$ is the indicator function of the set $A$. Clearly, $f_\delta \to f$ a.s. when $\delta \to 0$ as $\xi_n > 0$ a.s. Let $\delta' < \delta$. We will need a variation of $f_\delta$ for which the points are normalized with respect to the sum of the points within $[\delta', 1]$

$$f_{\delta, \delta'}(\xi_1, \ldots, \xi_n; Q_n) := f_\delta(\xi_1 / N_{\delta'}, \ldots, \xi_n / N_{\delta'}; Q_n),$$

where

$$N_{\delta'} := \sum_{i \,:\, \xi_i \geq \delta'} \xi_i.$$

Note that $N_{\delta'} \to 1$ almost surely when $\delta' \to 0$. Therefore, by continuity of $f$ and the fact that

$$\lim_{\delta' \to 0} \chi_{\{\xi_n / N_{\delta'} \geq \delta\}} \to \chi_{\{\xi_n \geq \delta\}}$$

we have

$$\lim_{\delta \to 0} \lim_{\delta' \to 0} f_{\delta, \delta'}(\xi_1, \ldots, \xi_n; Q_n) = f(\xi_1, \ldots, \xi_n; Q_n) \qquad \text{a.s.}$$

In particular, by (4.1) and the dominated convergence theorem,

$$(4.3) \qquad \lim_{\delta \to 0} \lim_{\delta' \to 0} \lim_{r \to \infty} \mathbb{E}_r[f_{\delta, \delta'}(\Phi_r(\xi, Q))] = \mathbb{E}[f(\xi, Q)].$$

On the other hand, let $A_{N, \delta'}^c$ be the event that all evolved points in $[\delta', 1]$ come from the first $N$ before evolution (as defined in Lemma 4.1 with $\lambda = 1$ and $T = 1$). Let us write $\Phi_r(\xi, Q)|_N$ for the evolution of the first $N$ points of $(\xi, Q)$. In this notation, the function $f_{\delta,\delta'}(\Phi_r(\xi, Q))$ when restricted to the event $A^c_{N,\delta'}$ only depends on $\Phi_r(\xi, Q)|_N$. Hence, by writing $1 = \chi_{A^c_{N,\delta'}} + \chi_{A_{N,\delta'}}$

$$
\begin{aligned}
(4.4) \qquad \mathbb{E}[f_{\delta,\delta'}(\Phi_r(\xi, Q))] &= \mathbb{E}[f_{\delta,\delta'}(\Phi_r(\xi, Q)|_N)\chi_{A^c_{N,\delta'}}] \\
&\quad + \mathbb{E}[f_{\delta,\delta'}(\Phi_r(\xi, Q))\chi_{A_{N,\delta}}].
\end{aligned}
$$

Clearly, the restriction to $N$ points of $Q^{*r}$ tends to the identity matrix. Recall that convergence in distribution for Gaussian vectors is equivalent to convergence of the covariance matrix thus

$$
(4.5) \qquad \lim_{r \to \infty} \mathbb{E}[f_{\delta,\delta'}(\Phi_r(\xi, Q)|_N)] = \mathbb{E}[f_{\delta,\delta'}(\Phi_\infty(\xi, Q)|_N)].
$$

We can take $N$ arbitrary large in (4.4) so that the probability of $A_{N,\delta}$ is small uniformly in $r$ by Lemma 4.1 thus

$$
\lim_{r \to \infty} \mathbb{E}[f_{\delta,\delta'}(\Phi_r(\xi, Q))] = \mathbb{E}[f_{\delta,\delta'}(\Phi_\infty(\xi, Q))].
$$

Equation 4.2 is obtained from (4.3) by taking the limit $\delta, \delta' \to 0$ of the above. $\square$

We stress that although robust quasi-stationarity implies quasi-stationarity under free evolution, it is generally not true that structures that are ergodic under correlated evolution will be ergodic under the free evolution. Indeed, recall that for the free evolution $\Phi_\infty$ the measure $\mu$ was $\Phi_\infty$-invariant and therefore deterministic if $(\xi, Q)$ is ergodic. We will see that this measure is never invariant under $\Phi_r$. Understanding how this measure transforms under correlated evolution is the object of the next section.

4.2. *The directing random overlap structure.* A combination of Theorems 1.10 and 1.9 guarantees the existence of a directing measure on a Hilbert space $\mathcal{H}$ in the case where $(\xi, Q)$ is robustly quasi-stationary. Without loss of generality, we can assume that $\mathcal{H}$ is a Gaussian Hilbert space. This perspective yields a representation of the Gaussian field $\kappa$ with covariance $Q$

$$
\kappa_i = \phi_i + \kappa_i^f,
$$

where $(\kappa_i^f, i \in \mathbb{N})$ are independent centered Gaussian variables of variance $1 - \|\phi_i\|^2$ depending on $i$. In the case where $\kappa$ has covariance $Q^{*r}$, $r \in \mathbb{N}$, we associate to each $\phi_i$ a vector $\phi_i^{(r)}$ in such a way that

$$
(4.6) \qquad (\phi_i^{(r)}, \phi_j^{(r)}) = (\phi_i, \phi_j)^r.
$$

This map is unique up to isometry of $\mathcal{H}$. We then have the representation $\kappa_i = \phi_i^{(r)} + \kappa_i^f$ where $(\kappa_i^f, i \in \mathbb{N})$ are independent centered Gaussian variables of variance $1 - \|\phi_i\|^r$.

Precisely, the combination of the two results yields:

**Theorem 4.2.** *Let $(\xi, Q)$ be a ROSt that is robustly quasi-stationary and ergodic for a function $\psi$ satisfying Assumption 1.3. The following hold:*

1. *$\xi$ is a $PD(x, 0)$ independent of $Q$ for some $x \in (0, 1)$;*
2. *$Q$ is directed by $\mu$, a random probability measure on $\mathcal{H}$ whose law, up to isometry, satisfies the equation*

$$(4.7) \qquad \frac{\mu(d\phi) e^{\psi_{x, \|\phi\|^2}(\phi)}}{\int_{\mathcal{H}} \mu(d\phi) e^{\psi_{x, \|\phi\|^2}(\phi)}} \stackrel{\mathcal{D}}{=} \mu(d\phi),$$

   *where*

$$(4.8) \qquad e^{\psi_{x, \rho}(y)} := \int_{\mathbb{R}} \frac{e^{-z^2/2}}{\sqrt{2\pi}} e^{x\psi(y + z\sqrt{1-\rho})} \, dz.$$

Proof. Theorem 1.9 together with Theorem 1.10 ensures that the law of $\xi$ is a superposition of $PD(x, 0)$ as well as the existence of the directing measure $\mu$.

From the remark preceding the theorem, the correlated evolution $\Phi_r$ reduces to the evolution of a $(\mathcal{H}, \mu)$-marking of a $PD(x, 0)$ variable for the choice $W_\phi(\kappa_i^f) := e^{\psi(\phi^{(r)} + \kappa_i^f)}$ (conditionally on $\mu$ and on the randomness of the Gaussian Hilbert space $\mathcal{H}$). Recall from Corollary 2.8 that a marked $PD(x, 0)$ under such a stochastic map is again a marked $PD(x, 0)$ with modified mark probability. In particular, the $x$ parameter is $\Phi_r$-invariant. Since $(\xi, Q)$ is ergodic and $PD(x, 0)$ with distinct parameters are mutually singular, we conclude that $\xi$ is simply a $PD(x, 0)$ for some $x \in (0, 1)$ independently of $\mu$ and the randomness of $\mathcal{H}$. The first claim is proven. Equation (4.7) is a consequence of quasi-stationarity and the form of $\bar{\mu}$ in Corollary 2.8. $\square$

In the next section, we solve the fixed point (4.7) in the particular case where $Q$ is indecomposable with $S_Q$ finite. It is proven in Theorem A.6 of Appendix A that under these assumptions the directing measure is discrete, of the form

$$\sum_{l \in \mathcal{L}} \tilde{\xi}_l \delta_{\phi_l},$$

where the index set $\mathcal{L}$ is countable, $\tilde{\xi}_1 \geq \tilde{\xi}_2 \geq \cdots > 0$ and the vectors $(\phi_l, l \in \mathcal{L}) \subset$ have square norm $q_{\max} := \max S_Q$.

We construct a ROSt that is a representation of $\mu$ up to isometry. Plainly $\tilde{\xi} := (\tilde{\xi}_l, l \in \mathcal{L})$ is a random mass-partition. The natural choice for the overlaps of $\tilde{\xi}$ is

$$\tilde{q}_{kl} := \frac{1}{q_{max}}(\phi_k, \phi_l),$$

where $q_{max} := \max S_Q$. Note that $q_{max} = 0$ if and only if $S_Q = \{0\}$ as the vectors $\phi_i$ must then have zero norm. The factor $1/q_{max}$ is added so that the diagonal elements are normalized to 1. Note that by definition

$$S_{\tilde{Q}} = \frac{1}{q_{max}}(S_Q \setminus \{q_{max}\}),$$

so that $|S_{\tilde{Q}}| = |S_Q| - 1$. We say that $(\tilde{\xi}, \tilde{Q})$ is the *directing ROSt of* $(\xi, Q)$.

In the discrete case, the law of $(\xi, Q)$ is equivalent to the law of a $(\mathcal{C}, \tilde{\xi})$-marking of $\tilde{\xi}$ with $\mathcal{C} := (\phi_l, l \in \mathcal{L}) \subset \mathcal{H}$. We write

(4.9)                         $\xi^{\mathcal{C}} := ((\xi_i, \phi_{i^*}), i \in \mathbb{N})$

for this marking where $(i^*, i \in \mathbb{N})$ are i.i.d. $\tilde{\xi}$-distributed. The striking fact is that the evolution of $\mu$ in (4.7) directly translates into a correlated evolution of $(\tilde{\xi}, \tilde{Q})$

(4.10)                $\left( \left( \frac{\tilde{\xi}_l e^{\psi_{x,q^r_{max}}(\sqrt{q^r_{max}}\phi_l^{(r)})}}{\sum_{l' \in \mathcal{L}} \tilde{\xi}_{l'} e^{\psi_{x,q^r_{max}}(\sqrt{q^r_{max}}\phi_{l'}^{(r)})}}, l \in \mathcal{L} \right)_\downarrow, \tilde{\pi}\tilde{Q}\tilde{\pi}^{-1} \right),$

where $\tilde{\pi}$ is obtained from the decreasing rearrangement of the evolved probabilities and the function $\psi_{x,\rho}(\sqrt{\rho}\cdot)$ is defined in (4.8). It is easily checked that this function satisfies Assumption 1.3.

By the above remarks, Theorem 4.2 becomes the following.

PROPOSITION 4.3. *Let $(\xi, Q)$ be a ROSt that is robustly quasi-stationary and ergodic for some function satisfying Assumption 1.3. If it is indecomposable and with finite state space, that is $S_Q = \{q_1, \ldots, q_k\}$ for some $-1 < q_1 < \cdots < q_k < 1$, with $k > 1$, then the following hold:*

1. *$\xi$ is a $PD(x, 0)$ independently of $Q$;*
2. *there exists a directing ROSt $(\tilde{\xi}, \tilde{Q})$ that is indecomposable with $S_{\tilde{Q}} = \{\frac{q_1}{q_k}, \ldots, \frac{q_{k-1}}{q_k}\}$ such that*

$$q_{ij} \overset{\mathcal{D}}{=} q_k \tilde{q}_{i^* j^*},$$

   *where $i^*, i \in \mathbb{N}$, are $\tilde{\xi}$-distributed random indices;*
3. *for all $r$ for which $(\xi, Q)$ is quasi-stationary, $(\tilde{\xi}, \tilde{Q})$ is quasi-stationary (and ergodic) under the correlated evolution $\Phi_{\psi_{x,q^r_k}}(\sqrt{q^r_k}\phi_l^{(r)})$ where $(\phi_l^{(r)}, l \in \mathcal{L})$ is centered Gaussian with covariance $\tilde{Q}^{*r}$ and $\psi_{x,\rho}$ is defined in (4.8).*

PROOF.  The law of $(\tilde{\xi}, \tilde{Q})$ must be invariant under the above evolutions for the same $r$ as $(\xi, Q)$ by robust quasi-stationarity. Moreover, $(\tilde{\xi}, \tilde{Q})$ ought to be ergodic as $(\tilde{\xi}, \tilde{Q})$ is $Q$-measurable. $\tilde{Q}$ is indecomposable as for all $i \in \mathbb{N}$

$$q_k S_{\tilde{Q}}(i^*) = S_Q(i) \setminus \{q_k\} = S_Q \setminus \{q_k\}. \qquad \square$$

4.3. *The hierarchical structure for linear* $\psi$.  We now use an inductive argument based on Proposition 4.3 to single out the RPC's in terms of quasi-stationarity properties for linear $\psi$. The induction is possible if one assumes robust quasi-stationarity under the evolution $\psi(\kappa) = \lambda\kappa$, where $\lambda = \lambda_r$ can depend on the power of the covariance matrix.

THEOREM 4.4.  *Let* $(\xi, Q)$ *be an indecomposable ROSt with* $|S_Q| = k$. *If* $(\xi, Q)$ *is robustly quasi-stationary, and ergodic under a sequence of evolutions with* $\psi(\kappa) = \lambda_r \kappa$, *for some uniformly bounded collection* $\{\lambda_r\}$, *then* $(\xi, Q)$ *is a* $k$-*level RPC*.

The strategic step is to prove that the directing ROSt $(\tilde{\xi}, \tilde{Q})$ inherits the same properties as $(\xi, Q)$. An induction on the cardinality of $S_Q$ completes the proof. An obstacle to this is the fact that the function for which $(\tilde{\xi}, \tilde{Q})$ is quasi-stationary is different for each field with covariance $\tilde{Q}^{*r}$ as seen from Proposition 4.3. In other words, $(\tilde{\xi}, \tilde{Q})$ is not robustly quasi-stationary for a fixed function $\psi$. However, this problem is circumvented for $\psi(\kappa) = \lambda\kappa$. In this case, it is easily checked from the definition of $\psi_{x,\rho}$ that the evolution is still linear since

$$(4.11) \qquad \psi_{x,\rho}(\sqrt{\rho}y) = x\sqrt{\rho}\lambda y + C.$$

The constant $C := x^2 \lambda^2 (1 - \rho)/2$ is irrelevant for the evolution because of the normalization. If $(\xi, Q)$ is assumed quasi-stationary for all $\lambda$, then so is $(\tilde{\xi}, \tilde{Q})$ since $x\sqrt{\rho}\lambda$ is a scaling of $\lambda$. A similar induction is still true when quasi-stationarity holds for a single $\lambda$ using several steps of the evolution.

LEMMA 4.5.  *Let* $(\xi, Q)$ *be a ROSt with a finite* $S_Q = \{q_1, \ldots, q_k\}$, *which:*

1. *is robustly quasi-stationary and ergodic under* $\Phi_r$ *with* $\psi(\kappa) = \lambda_r \kappa$, $\lambda_r$ *uniformly bounded,*
2. *is indecomposable,*
3. *has* $S_Q \subset (-1, 1)$.

*Then its directing ROSt* $(\tilde{\xi}, \tilde{Q})$, *with* $S_{\tilde{Q}} = \{\frac{q_1}{q_k}, \ldots, \frac{q_{k-1}}{q_k}\}$, *also satisfies the above three conditions.*

PROOF.  First we show that such ROSt's are quasi-stationary under the free evolution, that is, in the linear case the conclusion of Theorem 1.10 remains true even when $\lambda$ is allowed to depend on $r$. Indeed, consider $T$ time steps of the evolution with $T = T_r$. By the additive property of Gaussian variables, if $(\xi, Q)$ is quasi-stationary under $\psi(\kappa) = \lambda\kappa$, then it is quasi-stationary under $\lambda\sqrt{T}\kappa$ for all $T \in \mathbb{N}$. Since $\lambda_r$ is uniformly bounded, we can build a sequence $\lambda_r\sqrt{T_r}$ which converges to a nonzero limit as $r \to \infty$. The increments $\lambda_r\sqrt{T_r}\kappa_i$, where $\kappa$ has covariance $Q^{*r}$, then converge weakly to i.i.d. nonzero Gaussians. Moreover, the estimate of Lemma 4.1 is also uniform in $r$ since $\lim_{r\to\infty} g(\lambda_r)^{T_r}$ is finite. Quasi-stationarity under the free evolution then follows as in the proof of Theorem 1.10. In particular, following the same line of argument as in the previous section, we see that the conclusion of Proposition 4.3 remains valid when $\lambda$ depends on $r$ in the linear case. This proves that the directing ROSt satisfies properties 1 and 2.

It remains to show that $S_{\tilde{Q}} \subset (-1, 1)$. By definition of $\tilde{Q}$, it plainly follows that $\tilde{q}_{ij} \neq 1$ when $i \neq j$. Moreover, for a fixed $i$, there exists at most one $j$ for which $\tilde{q}_{ij} = -1$, so $\kappa_i(t) = -\kappa_j(t)$ for all $t$. It is shown in [6] that the time-average $v_i$ of the past increment must be the same for all points and is nonzero unless there exists only one point. This contradicts the existence of a pair $\tilde{q}_{ij} = -1$ since in that case $v_i = -v_j$. The lemma is proven.  □

We are now ready to prove Theorem 4.4. The proof is based on the description of RPC's given in Theorem 2.9.

PROOF OF THEOREM 4.4.   The proof is by induction on $|S_Q|$. Recall that a 1-level RPC is simply a pair $(\xi, Q)$ where $\xi$ is $PD(x, 0)$ and $Q$ is a constant matrix with all nondiagonal entries equal to $q$ for some $0 \leq q < 1$. Let $(\xi, Q)$ be as in the assumption of the theorem with $S_Q = \{q\}$ for $-1 < q < 1$. Note that the Gaussian field $\kappa$ with covariance $Q$ is exchangeable. In particular, $q$ is nonnegative. We can write

$$\kappa_i = \kappa_c + \kappa_i^f,$$

where $\kappa_c$ is a centered Gaussian with variance $q$ common to every point and $(\kappa_i^f, i \in \mathbb{N})$ are i.i.d. centered Gaussians of variance $1 - q$. As $\psi$ is linear, $\kappa_c$ is irrelevant in the correlated evolution as it cancels in the normalization being common to every point. We conclude that $(\xi, Q)$ is quasi-stationary and ergodic under the free evolution. By Theorem 3.1, $\xi$ is $PD(x, 0)$ and we conclude that $(\xi, Q)$ is a 1-level RPC.

Let $|S_Q| = k > 1$ and assume the result holds in the case $S_Q$ takes $k - 1$ values. By Proposition 4.3, $\xi$ has Poisson–Dirichlet distribution for some parameter $x_k \in (0, 1)$ independently of $Q$. Moreover the law of $Q$ is determined by a directing ROSt $(\tilde{\xi}, \tilde{Q})$ with $S_{\tilde{Q}} = \{q_1/q_k, \ldots, q_{k-1}/q_k\}$. In view of Theorem 2.9 on the law of the cascade, it remains to show that $(\tilde{\xi}, \tilde{Q})$ is a $(k-1)$-level RPC with parameters $0 < x_1/x_k < \cdots < x_{k-1}/x_k < 1$ for some $0 < x_1 < \cdots < x_{k-1} < 1$ and that $0 \le q_1/q_k < \cdots < q_{k-1}/q_k < 1$. We know from Lemma 4.5 that $(\tilde{\xi}, \tilde{Q})$ satisfies the assumptions of the theorem. By the induction hypothesis, it is a $(k-1)$-level RPC for some parameters $0 < x'_1 < \cdots < x'_{k-1} < 1$ and $0 \le q_1/q_k, \ldots, q_{k-1}/q_k < 1$. The proof is completed by defining $x_l$ as $x_l := x_k x'_l$ for $l = 1, \ldots, k-1$. The theorem is proven.  □

4.4. *The hierarchical structure for general $\psi$.*  A slight adaptation of the proof of the linear case yields the main result of this work, as it is stated in Theorem 1.8. The proof is based on the evolution of the systems after $T$ independent steps of the evolution taking $T$ to infinity. An application of the central limit theorem then permits to reduce the evolution to the linear case.

Throughout this section, we will assume without loss of generality that

$$\int_{\mathbb{R}} \frac{e^{-z^2/2}}{\sqrt{2\pi}} \psi(z) \, dz = 0$$

and

$$\int_{\mathbb{R}} \frac{e^{-z^2/2}}{\sqrt{2\pi}} \psi^2(z) \, dz = 1.$$

Note that we can always substract a constant to $\psi$ as the evolution is normalized and divide the function by its second moment as quasi-stationarity holds for multiples of $\psi$ in a neighborhood of 0.

We define the covariance function $C_\psi$

$$C_\psi(q) := \mathbb{E}[\psi(X)\psi(Y)],$$

where $X$ and $Y$ are standard Gaussians with covariance $q$. It is easy to see that $C_\psi$ is a continuous function on $[-1, 1]$ and that $\lim_{q \to 0} C_\psi(q) = 0$. We can actually say more when $\psi$ is a smooth function by a simple Gaussian differentiation.

LEMMA 4.6.   *If $\psi : \mathbb{R} \to \mathbb{R}$ is in $C^2(\mathbb{R})$ with bounded derivatives, then*

$$\frac{d}{dq} C_\psi(q) = C_{\psi'}(q).$$

*In particular, when $\psi$ is not a constant, the function $q \mapsto C_\psi(q)$ is strictly monotone in a neighborhood of 0.*

Let $\kappa = (\kappa_i, i \in \mathbb{N})$ be the centered Gaussian field with covariance matrix $Q^{*r}$. The covariance matrix $\hat{Q}(r)$ of the field $(\psi(\kappa_i), i \in \mathbb{N})$

(4.12) $$\hat{q}_{ij}(r) := C_\psi(q_{ij}^r).$$

It turns out that the ROSt $(\xi, \hat{Q}(r))$ inherits the quasi-stationarity property from $(\xi, Q)$.

LEMMA 4.7.   *Let $\psi$ be a function satisfying Assumption 1.3. If $(\xi, Q)$ is quasi-stationary under $\Phi_r$ for every multiple $\lambda\psi$, $\lambda$ in a neighborhood of 0, then $(\xi, \hat{Q}(r))$ is quasi-stationary for the linear functions $\lambda\hat{\kappa}$. Moreover, if $(\xi, Q)$ is indecomposable, then so is $(\xi, \hat{Q}(r))$.*

PROOF.   The proof is identical to the proof of Theorem 1.10 in Section 4.1. It suffices to replace $\Phi_r$ by

$$\Phi_T(\xi, Q) := \Phi_{\lambda/\sqrt{T}\psi(\kappa(T-1))} \circ \cdots \circ \Phi_{\lambda/\sqrt{T}\psi(\kappa(0))}.$$

The requirement that quasi-stationarity holds for $\lambda$ in a neighborhood of 0 ensures that $(\xi, Q)$ is quasi-stationary after the scaling $\lambda \to \lambda/\sqrt{T}$. The probability of the event $A_{N,T,\lambda/\sqrt{T},\delta}$ of Lemma 4.1 is uniformly small in $T$ when $N$ is large since $\lim_{T\to\infty} g(\lambda/\sqrt{T})^T < \infty$. The proof of Theorem 1.10 thus applies. The equivalent of (4.5) is obtained through the $N$-dimensional central limit theorem

$$\lim_{T\to\infty} \mathbb{E}[f_{\delta,\delta'}(\Phi_T(\xi,Q)|_N)] = \mathbb{E}[f_{\delta,\delta'}(\Phi_{\lambda\hat{\kappa}}(\xi,\hat{Q})|_N)],$$

where $\hat{\kappa}$ is a centered Gaussian field with covariance $\hat{Q}(r)$. Indecomposability is clear from the fact that

$$S_{\hat{Q}(r)}(i) = C_\psi(S_{Q^{*r}}(i)) = C_\psi(S_{Q^{*r}}) = S_{\hat{Q}(r)}. \qquad \square$$

We now have all the tools to prove Theorem 1.8.

PROOF OF THEOREM 1.8.   By Lemma 4.7, $(\xi, \hat{Q}(r))$ is robustly quasi-stationary under evolutions corresponding to the linear functions $\lambda\hat{\kappa}$ for $\lambda$ in a neighborhood of 0. Due to the $Q$-factorization, we can assume that $(\xi, \hat{Q}(r))$ is indecomposable. As discussed above Lemma 4.6,

$$\lim_{r\to\infty} \hat{q}(r)_{ij} = \delta_{ij}.$$

Hence the proof of Theorem 4.4 holds *verbatim* with $Q^{*r}$ replaced by $\hat{Q}(r)$. We conclude that $(\xi, \hat{Q})$ is a $k$-level RPC. It remains to prove that so is $(\xi, Q)$. Recall that

(4.13)                     $$\hat{q}(r)_{ij} = C_\psi(q_{ij}^r).$$

By Lemma 4.6, $C_\psi$ is strictly monotone in a neighborhood of 0, say $V$. In particular, it is invertible on $V$ and $C_\psi^{-1} : V \to [0, 1]$ is also strictly monotone. By choosing $r$ large enough so that $(\max S_Q)^r$ belongs to $V$, we can invert (4.13)

$$q_{ij}^r = C_\psi^{-1}(\hat{q}(r)_{ij}).$$

By Proposition 2.4, we conclude that $(\xi, Q^{*r})$ is a $k$-level RPC and consequently, that $(\xi, Q)$ is. The theorem is proven.   $\square$

## APPENDIX A: WEAKLY EXCHANGEABLE MATRICES

The goal of this section is to prove that the directing measure of a weakly exchangeable overlap matrix must be discrete, under the assumption of indecomposability and finiteness of the state space. This is done by defining a random partition of $\mathbb{N}$ from the random overlap matrix. Proof of the discreteness of the measure is then reduced to showing that the blocks of the partition are of nonzero density.

**A.1. Weakly exchangeable covariance matrices.** A random $\mathbb{N} \times \mathbb{N}$-matrix $M$ is said to be *weakly exchangeable* if

$$\tau M \tau^{-1} \overset{\mathcal{D}}{=} M$$

for any permutation $\tau$ of finite elements of $\mathbb{N}$. We say that a matrix is a *covariance matrix* in the case it is a real positive semi-definite symmetric matrix. A variation of de Finetti's theorem due to Dovbysh and Sudakov characterizes the random covariance matrices that are weakly exchangeable (see also a proof by Hestir [13]).

THEOREM A.1 [10]. *Let $C = \{c_{ij}\}$ be a random $\mathbb{N} \times \mathbb{N}$ covariance matrix that is weakly exchangeable and for which $\mathbb{E}[c_{11}] < \infty$. Then there exists a random probability measure $\mu$ on $\mathcal{H} \times [0, \infty)$, with $\mathcal{H}$ some separable Hillbert space, such that the law of $C$ is given by*

$$c_{ij} \overset{\mathcal{D}}{=} (\phi_i, \phi_j)_{\mathcal{H}} + \delta_{ij} a_i,$$

*where $\{(\phi_i, a_i)\}_{i \in \mathbb{N}}$ is an i.i.d. $\mu$-distributed sequence.*

We say that $\mu$ is the *directing measure of $C$*. The above result is sharp since any probability measure on $\mathcal{H} \times [0, \infty)$ can be used to construct a weakly exchangeable covariance matrix.

In the notation of the above theorem, a weakly exchangeable random overlap matrix $Q$ (i.e., a covariance matrix with $q_{ii} = 1$ for all $i$) is such that $\|\phi_i\|^2 + a_i = 1$ and so the diagonal part $a_i$ is determined by $\phi_i$. In particular, the directing $\mu$ can be simply seen as a probability measure on $\mathcal{H}$.

**A.2. Partitions of $\mathbb{N}$ and weak exchangeability.** There is a natural exchangeable random partition of $\mathbb{N}$ that can be constructed from a weakly exchangeable overlap matrix $Q$.

Recall that a *partition $\Gamma$ of $\mathbb{N}$* is a countable collections of disjoint subsets $\Gamma := (\Gamma_j, j \in \mathbb{N})$, also called *blocks*, for which

$$\bigcup_{j \in \mathbb{N}} \Gamma_j = \mathbb{N}.$$

In this notation, we set $\Gamma_j = \varnothing$ for $j$ large enough in the case the partition has finite number of blocks. Plainly, a partition $\Gamma$ is determined by the equivalence relation

$$i \sim j \iff i \text{ and } j \text{ belong to the same block of } \Gamma.$$

We say that a block $\gamma \in \Gamma$ has a *density* in case the following limit exists

$$|\gamma| := \lim_{n \to \infty} \frac{1}{n} \#(\gamma \cap \{1, \dots, n\}).$$

Random partitions of $\mathbb{N}$ are random variables on the space of partitions of $\mathbb{N}$, which is equipped with the $\sigma$-algebra generated by the sets $\{\Gamma : i \sim j\}$ for $i, j \in \mathbb{N}$. The laws are constructed from consistent distributions on the space of partitions of $\{1, \dots, n\}$. A random partition of $\mathbb{N}$ is said to be *exchangeable* if and only if for all $n \in \mathbb{N}$, its law restricted to $\{1, \dots, n\}$ is the same after any permutation of the $n$ elements.

An example of exchangeable random partitions is given by the so-called *paintbox based on a mass-partition $s$* [7]. The partition of $\mathbb{N}$ is constructed by pairing the weights of $s$ to disjoint intervals of $[0, 1]$ in such a way that the length of the interval $I_i$ equals the weight $s_i$. Taking $(U_i, i \in \mathbb{N})$ i.i.d. uniform random variables on $[0, 1]$, the partition is defined by the equivalence relation

$$i \sim j \iff U_i \text{ and } U_j \text{ belongs to the same interval of the partition.}$$

It is easy to see that each block of a paintbox has a density corresponding to the corresponding weight of $s$. Moreover, the densities are strictly positive if and only if $s$ is proper, that is $\sum_i s_i = 1$. Exchangeable random partitions of $\mathbb{N}$ are characterized in the following theorem of Kingman.

THEOREM A.2 [7]. *Let $\Gamma$ be an exchangeable random partition of $\mathbb{N}$. Then its law is a linear superposition of paintboxes. Equivalently, $\Gamma$ is paintbox based on a random mass-partition.*

Following Kingman's idea, we construct an exchangeable random partition of $\mathbb{N}$ from a weakly exchangeable overlap matrix $Q$ as follows. In view of Theorem A.1, we define the equivalence relation

$$i \sim j \iff \phi_i = \phi_j,$$

where $(\phi_i, i \in \mathbb{N})$ is the collection of random elements of $\mathcal{H}$ whose existence is asserted by the theorem. We write $\Gamma_Q$ for the random partition of $\mathbb{N}$ induced by this equivalence relation. As $(\phi_i, i \in \mathbb{N})$ is exchangeable, so is this partition. A straightforward application of Theorem A.2 shows that the law of $\Gamma_Q$ is a superposition of paintboxes. However, there is no reason a priori for the random mass-partition on which $\Gamma_Q$ is based to be proper. It would be so if the blocks of $\Gamma_Q$ have positive density a.s. We introduce two sufficient conditions on the overlap matrix for this to hold.

ASSUMPTION A.3.  The set of values taken by the nondiagonal entries of the $i$th row, $S_Q(i) := \{q_{ij} : i > j\}$, satisfy the following:

1. $|S_Q(i)| < \infty$ for all $i \in \mathbb{N}$ $Q$-a.s.,
2. $S_Q(i) = S_Q(j)$ for all $i, j \in \mathbb{N}$ $Q$-a.s., in which case we write $S_Q$ for the common set of values.

In the simple case where $S_Q = \{q\}$, the overlap matrix must be of the form $q_{ij} = q(1 - \delta_{ij}) + \delta_{ij}$. Therefore the restriction of the measure $\mu$ to the space $\mathcal{H}$ is entirely supported on a single vector of square norm $q$. In particular, $q$ must be nonnegative. The next results generalize the above reasoning to the case where $S_Q$ is finite.

LEMMA A.4.  *Let $Q$ be a weakly exchangeable overlap matrix, $\Gamma_Q$ its associated random partition of $\mathbb{N}$, and $(\phi_i, i \in \mathbb{N})$ the random sequence determining $Q$ in Theorem A.1. If $Q$ satisfies Assumption A.3(1), then*

$$\|\phi_i\|^2 \in S_Q(i),$$

*and in particular $\max S_Q(i)$ is nonnegative. Moreover, almost surely:*

$$\text{(A.1)} \qquad \lim_{n \to \infty} \frac{1}{n} \#\{i + 1 \le j \le i + n : q_{ij} = \|\phi_i\|^2\} > 0.$$

PROOF.  Consider for fixed $i$ the sequence $((\phi_i, \phi_j), j > i))$. The first claim will be proven if for all $\delta > 0$, there exists almost surely a $j > i$ such that

$$\text{(A.2)} \qquad |(\phi_i, \phi_j) - \|\phi_i\|^2| < \delta.$$

Indeed, this would show that the nondiagonal entries come arbitrarily close to the square norm and the result follows by the finiteness of $S_Q$. Equation (A.2) will be established if for all $i \in \mathbb{N}$ and $\delta > 0$

$$\text{(A.3)} \qquad \mu(\{\phi \in \mathcal{H} : |(\phi, \phi_i) - \|\phi_i\|^2| < \delta\}) > 0$$

as we will have

$$\bigotimes_{j > i} \mu\left(\bigcap_{j > i} \{|(\phi_i, \phi_j) - \|\phi_i\|^2| \ge \delta\}\right) = \prod_{j > i} \mu(\{\phi_j : |(\phi_i, \phi_j) - \|\phi_i\|^2| \ge \delta\})$$

$$= \lim_{N \to \infty} \mu(\{\phi : |(\phi_i, \phi) - \|\phi_i\|^2| \ge \delta\})^N$$

$$= 0.$$

Let $(\phi_i, i \in \mathbb{N})$ be a sequence of i.i.d. $\mu$-distributed. Note that for all $\varepsilon > 0$, we must have

$$\text{(A.4)} \qquad \mu(\{\phi \in \mathcal{H} : \|\phi - \phi_i\| < \varepsilon\}) > 0 \qquad (\phi_i, i \in \mathbb{N})\text{-a.s.}$$

for all $i \in \mathbb{N}$. On the other hand, for all $\varepsilon > 0$ there exists a $\delta$ such that

$$\{\phi \in \mathcal{H} : \|\phi - \phi'\| < \varepsilon\} \subseteq \{\phi \in \mathcal{H} : |(\phi, \phi') - \|\phi'\|^2 < \delta\}.$$

Equation (A.3) is obtained from the above and (A.4). Equation (A.1) is obvious from the fact that $(\phi_i, \phi_j)$ takes a finite number of values and the sequence $((\phi_i, \phi_j), j > i)$ is i.i.d. $\quad\square$

LEMMA A.5. *Let $Q$ be as in Lemma A.4. If Assumption A.3*(ii) *is also satisfied, then the following hold:*

1. *The directing measure $\mu$ on $\mathcal{H}$ is almost surely supported on elements with uniform square norm $\max S_Q$*

$$\|\phi_i\|^2 = \max S_Q \qquad \text{for all } i \in \mathbb{N}.$$

   *In particular, $i$ and $j$ belong to the same block of $\Gamma_Q$ if and only if $q_{ij} = \max S_Q$.*
2. *Every block of $\Gamma_Q$ has positive density. In particular, $\Gamma_Q$ is a paintbox based on a proper random mass-partition.*

PROOF.    Fix $i \in \mathbb{N}$. From Cauchy–Schwarz, it is clear that for all $j \in \mathbb{N}$

$$q_{ij} \leq \|\phi_i\|\|\phi_j\| \leq \max S_Q$$

as $\|\phi_k\|^2 \in S_Q$ for all $k \in \mathbb{N}$ by Lemma A.4. Moreover, by Assumption A.3(ii) there must exist a $j > i$ such that $q_{ij} = \max S_Q$. Therefore from the above equation

$$\max S_Q = \max_j q_{ij} \leq \|\phi_i\|(\max S_Q)^{1/2} \leq \max S_Q.$$

We conclude that $\|\phi_i\|^2 = \max S_Q$. Furthermore, from

$$0 \leq \|\phi_i - \phi_j\|^2 = 2\max S_Q - 2q_{ij},$$

we deduce that $i$ and $j$ are in the same block of $\Gamma_Q$ if and only if $q_{ij} = \max S_Q$.

The last claim is obtained from (A.1) and the first assertion.    $\square$

The main result of this section gives a description of the directing measure $\mu$ under Assumption A.3: it must be supported on a countable collection of vectors.

THEOREM A.6.    *Let $Q$ be a weakly exchangeable overlap matrix satisfying Assumptions A.3. Then the directing measure $\mu$ of $Q$ on $\mathcal{H}$ is of the form*

$$\mu = \sum_{l \in \mathcal{L}} p_l \delta_{\phi_l}.$$

*The index set $\mathcal{L}$ is countable, $(\phi_l, l \in \mathcal{L})$ is a collection of elements of $\mathcal{H}$ of square norm $\max S_Q$, and $p := (p_l, l \in \mathcal{L})$ with $p_1 \geq p_2 \geq \cdots > 0$ are the corresponding probability weights.*

PROOF. We recall that $\Gamma_Q$ is based on a proper random mass-partition say $p$. Consider $(p_l, l \in \mathcal{L})$ ordered in a decreasing order of $p_i > 0$. By definition of the equivalence class of $\Gamma_Q$, for each class there is an identifying vector in the Hilbert space $\mathcal{H}$. Reindex these elements in decreasing order of probability (with multiplicities resolved through an ordering of $\mathcal{H}$). In this way we obtain the collection $(\phi_l, l \in \mathcal{L})$ on which $p$ is supported. □

**Acknowledgments.** We thank Jason Miller for valuable comments on the first version of the paper.

## REFERENCES

[1] AIZENMAN, M., SIMS, R. and STARR, S. L. (2003). An extended variational principle for the SK spin-glass model. *Phys. Rev. B* **68** 214403.

[2] AIZENMAN, M., SIMS, R. and STARR, S. L. (2007). Mean-field spin glass models from the cavity-ROSt perspective. In *Prospects in Mathematical Physics. Contemporary Mathematics* **437** 1–30. Amer. Math. Soc., Providence, RI. MR2354653

[3] ALDOUS, D. J. (1985). Exchangeability and related topics. In *École d'été de Probabilités de Saint–Flour, XIII—1983. Lecture Notes in Math.* **1117** 1–198. Springer, Berlin. MR883646

[4] ARGUIN, L.-P. (2007). Spin glass computations and Ruelle's probability cascades. *J. Stat. Phys.* **126** 951–976. MR2311892

[5] ARGUIN, L.-P. (2007). A dynamical characterization of Poisson–Dirichlet distributions. *Electron. Comm. Probab.* **12** 283–290 (electronic). MR2342707

[6] ARGUIN, L.-P. (2008). Competing particle systems and the Ghirlanda–Guerra identities. *Electron. J. Probab.* **13** 2101–2117. MR2461537

[7] BERTOIN, J. (2006). *Random Fragmentation and Coagulation Processes. Cambridge Studies in Advanced Mathematics* **102** 288. Cambridge Univ. Press, Cambridge. MR2253162

[8] BOLTHAUSEN, E. and SZNITMAN, A.-S. (1998). On Ruelle's probability cascades and an abstract cavity method. *Comm. Math. Phys.* **197** 247–276. MR1652734

[9] DERRIDA, B. (1985). A generalization of the random energy model which includes correlations between energies. *J. Phys. Lett.* **46** L401–L407.

[10] DOVBYSH, L. N. and SUDAKOV, V. N. (1982). Gram–de Finetti matrices. *J. Soviet. Math.* **24** 3047–3054. MR666087

[11] GUERRA, F. (2003). About the cavity fields in mean field spin glass models. Preprint. Available at arxiv:cond-mat/0307673.

[12] GUERRA, F. (2003). Broken replica symmetry bounds in the mean field spin glass model. *Comm. Math. Phys.* **233** 1–12. MR1957729

[13] HESTIR, K. (1989). A representation theorem applied to weakly exchangeable nonnegative definite arrays. *J. Math. Anal. Appl.* **142** 390–402. MR1014583

[14] HORN, R. A. and JOHNSON, C. R. (1985). *Matrix Analysis.* Cambridge Univ. Press, Cambridge. MR832183

[15] MÉZARD, M., PARISI, G. and VIRASORO, M. A. (1987). *Spin Glass Theory and Beyond. World Scientific Lecture Notes in Physics* **9**. World Scientific, Teaneck, NJ. MR1026102

[16] RUELLE, D. (1987). A mathematical reformulation of Derrida's REM and GREM. *Comm. Math. Phys.* **108** 225–239. MR875300

[17] RUZMAIKINA, A. and AIZENMAN, M. (2005). Characterization of invariant measures at the leading edge for competing particle systems. *Ann. Probab.* **33** 82–113. MR2118860

[18] SIMON, B. (1993). *The Statistical Mechanics of Lattice Gases, Vol. I.* Princeton Univ. Press, Princeton, NJ. MR1239893

[19] TALAGRAND, M. (2006). The Parisi formula. *Ann. Math. (2)* **163** 221–263. MR2195134

DEPARTMENT OF MATHEMATICS
PRINCETON UNIVERSITY
PRINCETON, NEW JERSEY 08544
USA
E-MAIL: larguin@math.princeton.edu

DEPARTMENTS OF MATHEMATICS AND PHYSICS
PRINCETON UNIVERSITY
PRINCETON, NEW JERSEY 08544
USA
E-MAIL: aizenman@math.princeton.edu